\documentclass[conference]{IEEEtran}
\IEEEoverridecommandlockouts

\usepackage{url}            
\usepackage{mathrsfs}
\usepackage{amssymb}
\usepackage{amsmath}
\usepackage{amsfonts}
\usepackage{amsthm}
\usepackage{algorithm}
\usepackage{algorithmic}
\usepackage{graphicx}
\usepackage{caption}
\usepackage[symbol]{footmisc}

\newtheorem{assumption}{Assumption}
\newtheorem{theorem}{Theorem}
\newtheorem{lemma}{Lemma}
\newtheorem{definition}{Definition}
\newtheorem{problem}{Problem}
\usepackage{hyperref}
\begin{document}
\title{Toward Model Parallelism for Deep Neural Network based on Gradient-free ADMM Framework}
\author{ Junxiang Wang\\Emory University\\ jwan936@emory.edu \and \IEEEauthorblockN{Zheng Chai,Yue Cheng}\IEEEauthorblockN{George Mason University}\IEEEauthorblockN{\IEEEauthorblockN{\{zchai2,yuecheng\}@gmu.edu}}\and Liang Zhao\IEEEauthorrefmark{1}\\
Emory University\\lzhao41@emory.edu\\ {*corresponding author}
}
\maketitle
\vspace{-0.7cm}
\begin{abstract}
  Alternating Direction Method of Multipliers (ADMM) has recently been proposed as a potential alternative optimizer to the Stochastic Gradient Descent(SGD) for deep learning problems. This is because ADMM can solve gradient vanishing and poor conditioning problems. Moreover, it has shown good scalability in many large-scale deep learning applications. However, there still lacks a parallel ADMM computational framework for deep neural networks because of layer dependency among variables. In this paper, we propose a novel parallel deep learning ADMM framework (pdADMM) to achieve layer parallelism: parameters in each layer of neural networks can be updated independently in parallel. The convergence of the proposed pdADMM to a critical point is theoretically proven under mild conditions. The convergence rate of the pdADMM is proven to be $o(1/k)$ where $k$ is the number of iterations. Extensive experiments on six benchmark datasets demonstrated that our proposed pdADMM can lead to more than 10 times speedup for training large-scale deep neural networks, and outperformed most of the comparison methods. Our code is available at: \textit{\url{https://github.com/xianggebenben/pdADMM}}.
\end{abstract}
\section{Introduction}
 \indent Due to wide applications and significant success in various applications, the training of deep neural network models has gained ever-increasing attention from the machine learning community. Gradient-based methods such as Stochastic Gradient Descent (SGD) and its variants have been considered as the state-of-the-art since the 1980s, mainly due to its superior performance. Despite the popularity, in recent years, the constant improvement of DNNs' performance is accompanied by a fast increase in models' complexity and size, which indicates a clear trend toward larger and deeper networks. Such a trend leads to severe challenges for large models to be fit into a single computing unit (e.g., GPU), and raises urgent demands for partitioning the model into different computing devices to parallelize training. However, the inherent bottleneck from backpropagation which prevents the gradients of different layers being calculated in parallel. This is because in backpropagation the gradient calculations of one layer tightly depend on and have to wait for the calculated results of all the previous layers, which prevents the gradients of different layers being calculated in parallel. 
 \\\indent 
 To work around the drawback from gradient-based methods, gradient-free methods have caught fast increasing attention in recent years, which aims to address the drawbacks such as strong dependency among layers, gradient vanishing (i.e. the error signal diminishes as the gradient is backpropagated), and poor conditioning (i.e. a small input can change the gradient dramatically). For example, Talyor et al. and Wang et al. presented an Alternating Direction Method of Multipliers (ADMM) algorithm to train neural network models \cite{taylor2016training,wang2019admm}. Moreover, extensive experiments have revealed that the ADMM outperformed most of SGD-related methods \cite{wang2019admm}. Amongst the gradient-free methods for deep learning optimization, ADMM-based methods are deemed to have great potential of parallelism of deep neural network training, due to its inherent nature, which can break an objective into multiple subproblems, each of which can be solved in parallel \cite{boyd2011distributed}.\\
 \indent Despite the potential, a parallel algorithm based on ADMM for deep neural network training has rarely been explored, developed, and evaluated until now, due to the layer dependency among subproblems of ADMM. Even though the ADMM reduces the layer dependency compared with SGD, one subproblem of ADMM is dependent on its previous subproblem. Therefore, existing ADMM-based optimizers still update parameters sequentially.\\
\indent To handle the difficulties of layer dependency,  in this paper we propose a novel parallel deep learning Alternating Direction Method of Multipliers (pdADMM) optimization framework to train large-scale neural networks. Our contributions in this paper include:
\begin{itemize}
\item We propose a novel reformulation of the feed-forward neural network problem, which splits a neural network into independent layer partitions and allows for ADMM to achieve model parallelism.
\item We present a model-parallelism version of the ADMM algorithm to train a feed-forward deep neural network. All parameters in each layer can be updated in parallel to speed up the training process significantly. All subproblems generated by the pdADMM algorithm are discussed in detail.
\item  We investigated the convergence properties of parallel ADMM in the common nonlinear activation functions such as the Rectified linear unit (Relu), and we prove that the pdADMM converges to a state-of-the-art critical point with a sublinear convergence rate $o(1/k)$.
\item We conduct extensive experiments on six benchmark datasets to show the massive speedup of the proposed pdADMM as well as its competitive performance with state-of-the-art optimizers.
\end{itemize}
\indent The organization of this paper is shown as follows: In Section \ref{sec:related work}, we summarize recent related research work to this paper. In Section \ref{sec:algorithm}, we formulate the novel pdADMM algorithm to train a feed-forward neural network. In Section \ref{sec:convergence},  the convergence guarantee of pdADMM to a critical point is provided. Extensive experiments on benchmark datasets to demonstrate the convergence, speedup and comparable performance of pdADMM are shown in Section \ref{sec:experiment}, and Section \ref{sec:conclusion} concludes this work.
\section{Related Work}
\label{sec:related work}
\indent \textbf{Distributed ADMM} ADMM is one of the commonly applied techniques in distributed optimization. Overall, the previous works on distributed ADMM can be classified into two categories:  synchronous problems and asynchronous problems. Synchronous problems usually require workers to optimize parameters in time before the master update the consensus variable, while asynchronous problems allow some workers to delay parameter updates. Most literature focused on the application of the distributed ADMM on synchronous problems. For example,  Mota et al. utilized the distributed ADMM for the congestion control problem \cite{mota2012distributed}; Makhdoumi and Ozdaglar studied the convergence properties of the distributed ADMM on the network communication problem. For more work, please refer to \cite{chang2016proximal,chang2014multi,shi2014linear,xu2017adaptive,zhu2016quantized}. On the other hand, a handful of papers investigated how asynchronous
problems can be addressed by distributed ADMM. For instance, Zhang et al., Wei et al., Chang et al and Hong proved the convergence of the distributed ADMM on asynchronous problems \cite{zhang2014asynchronous, wei2012distributed,chang2016asynchronous1,chang2016asynchronous2,hong2014distributed}. Kumar et al. discussed the application of the ADMM on multi-agent problems over heterogeneous networks \cite{kumar2016asynchronous}. However, there still lacks a general framework for ADMM to train deep neural networks in the distributed fashion.\\
\textbf{Convergence analysis of
 nonconvex ADMM:} Despite the outstanding performance of the nonconvex ADMM, its convergence theory is not well established due to the complexity of both coupled objectives and various (inequality and equality) constraints. Specifically, Magnusson et al. provided new convergence conditions of ADMM for a class of nonconvex structured optimization problems \cite{magnusson2015convergence}; Li and Pong investigated the properties of the nonconvex ADMM on the composite optimization problem \cite{li2015global}; Wang et al. presented mild convergence conditions of the nonconvex ADMM where the objective function can be coupled and nonsmooth \cite{wang2015global}; Hong et al. proved that the classic ADMM  converges to stationary points provided that the penalty parameter is sufficiently large \cite{hong2016convergence}; Wang et al. proved the convergence of multi-convex ADMM with inequality constraints \cite{wang2019multi}; Liu et al. proved the convergence properties of a parallel and linearized ADMM \cite{liu2019linearized}. Wang and Zhao studied the convergence conditions of the nonconvex ADMM in the nonlinearly constrained equality problems \cite{wang2017nonconvex}; Xie et al. proposed a deep-learning-based ADMM algorithm to study the constrained optimization problems \cite{xie2019differentiable}.  Wang et al. gave the first convergence proof of ADMM in the nonconvex deep learning problems \cite{wang2019admm,wang2019accelerated}. For more work, please refer to \cite{chartrand2013nonconvex,hajinezhad2015nonconvex,guo2017convergence,themelis2020douglas,wang2020tssm,wang2019opt}.\\
\indent \textbf{Distributed Deep Learning} With the increased volume of data and layers of neural networks, there is a need to design distributed systems to train a deep neural network for large-scale applications. Most recent papers have proposed gradient-based distributed systems to train neural networks: For example, Wen et al. proposed  Terngrad to accelerate distributed deep learning in data parallelism \cite{wen2017terngrad}; Sergeev et al. presented an open-source library Horovod to reduce communication overhead \cite{sergeev2018horovod}. Other systems include SINGA \cite{ooi2015singa} Mxnet\cite{chen2015mxnet}, TicTac \cite{hashemitictac} and Poseidon \cite{zhang2017poseidon}.\\
\textbf{Data and Model Parallelism} Data parallelism focuses on distributing data across different processors, which can be implemented in parallel. Scaling SGD is one of the most common ways to reach data parallelism \cite{zinkevich2010parallelized}. For example, the distributed architecture, Poseidon, is achieved by scaling SGD through overlapping communication and computation over networks.  The recently proposed ADMM \cite{taylor2016training,wang2019admm} is another way of data parallelism: each subproblem generated by ADMM can be solved in parallel. However, data parallelism suffers from the bottleneck of a neural network: for SGD, the gradient should be transmitted through all processors;  for ADMM, the parameters in one layer are subject to these in its previous layer. As a result, this leads to heavy communication cost and time delay. Model parallelism, however, can solve this challenge because model parallelism splits a neural network to many independent partitions. In this way, each partition can be optimized in parallel and hence reduce time delay. For instance, Parpas and Muir proposed a parallel-in-time method from the perspective of dynamic systems \cite{parpas2019predict}; Huo et al. introduced a feature replay algorithm to achieve model parallelism \cite{huo2018training}. Zhuang et al. broke layer dependency by introducing the delayed gradient \cite{zhuang2019fully}. However, to the best of our knowledge, there still lacks an exploration on how to achieve model parallelism via ADMM.
\section{pdADMM Algorithm}
\label{sec:algorithm}
\indent We propose the pdADMM algorithm in this section. Specifically, Section \ref{sec:problem setup} introduces the existing deep learning ADMM method, and reformulates the problem and presents the pdADMM algorithm in detail.
Section \ref{sec:subproblem} discusses all subproblems generated by pdADMM and the strategy to train a large-scale deep neural network via pdADMM.
\subsection{Background}
\label{sec:problem setup}
\begin{figure}
    \centering
    \includegraphics[width=\linewidth]{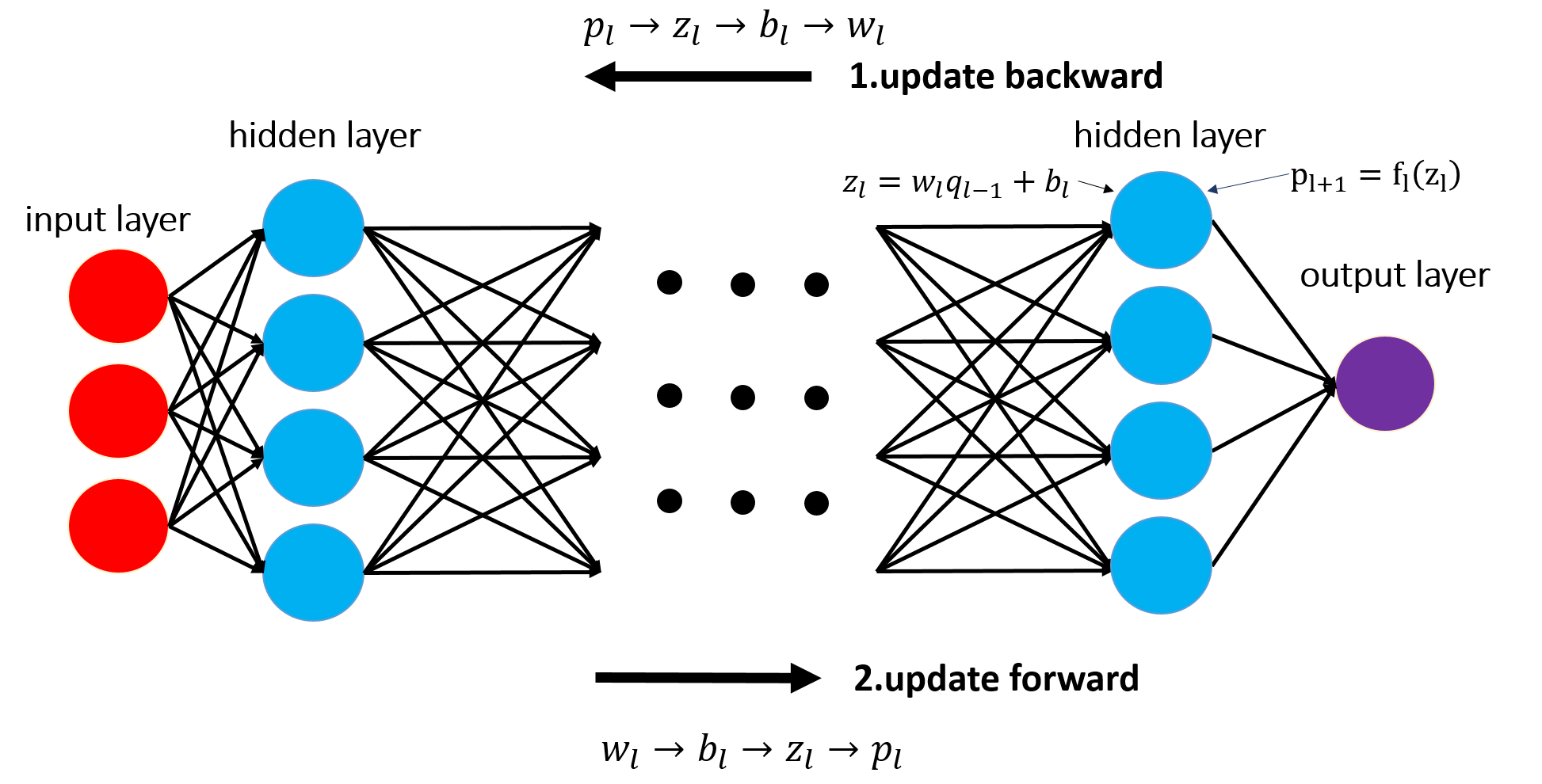}
    \caption{The overview of existing dlADMM algorithm: parameters are updated in a sequential fashion.}
    \label{fig:dlADMM}
\end{figure}
\begin{table}
\centering
 \begin{tabular}{cc}
 \hline
 Notations&Descriptions\\ \hline
 $L$& Number of layers.\\
 $W_l$& The weight matrix for the $l$-th layer.\\
 $b_l$& The intercept vector for the $l$-th layer.\\
 $z_l$& The auxiliary variable of the linear mapping for the $l$-th layer.\\
 $f_l(z_l)$& The nonlinear activation function for the $l$-th layer.\\
 $p_l$& The input for the $l$-th layer.\\
 $q_l$& The output for the $l$-th layer.\\
 $x$& The input matrix of the neural network.\\
 $y$& The predefined label vector.\\
 $R(z_L,y)$& The risk function for the $l$-th layer.\\
 $n_l$& The number of neurons for the $l$-th layer.\\
\hline
  \end{tabular}
  \captionof{table}{Important Notations}
   \label{tab:notation}

\end{table}

\indent In this section, we introduce the formulation of the feed-forward neural network training problem and an existing deep learning ADMM method. The important notations of this paper are detailed in Table \ref{tab:notation}.

The feed-forward neural network  is formulated as follows \cite{wang2019admm}:
\begin{problem}
\label{prob:problem 1}
\begin{align*}
     & \min\nolimits_{W_l,b_l,z_l,p_l} R(z_L;y) \\
     &s.t.\ z_l=W_lp_{l}+b_l, \   p_{l+1}=f_l(z_l)(l=1,\cdots,L-1)
\end{align*}
\end{problem}
where $p_1=x\in\mathbb{R}^{n_0}$ is the input of the deep neural network where $n_0$ is the number of feature dimensions, and $y$ is a predefined label vector. $p_l$ is the input for the $l$-th layer, also the output for the $(l-1)$-th layer. $R(z_L;y)$ is a risk function for the $L$-th layer, which is convex, continuous and proper. $z_l=W_lp_l+b_l$  and $p_{l+1}=f_l(z_l)$ are linear and nonlinear mappings for the $l$-th layer, respectively.\\
\indent Problem \ref{prob:problem 1} has been addressed by deep learning Alternating Direction Method of Multipliers (dlADMM) \cite{wang2019admm}. As shown in Figure \ref{fig:dlADMM}, the dlADMM algorithm updates parameters from the final layer, and  moves backward to the first layer, then updates parameters forward from the first layer to the final layer, in order to exchange information efficiently. However, parameters in one layer are dependent on its neighboring layers, and hence can not achieve parallelism. For example, the update of $p_{l+1}$ on the $l+1$-th layer needs to wait before $z_l$ on the $l$-th layer is updated.
\label{sec:pdADMM}
In order to address layer dependency, we relax Problem \ref{prob:problem 1} to Problem \ref{prob:problem 2} as follows:
\begin{problem}
\label{prob:problem 2}
\begin{align*}
    &\min\nolimits_{\textbf{p},\textbf{W},\textbf{b},\textbf{z},\textbf{q}} F(\textbf{p},\textbf{W},\textbf{b},\textbf{z},\textbf{q})=R(z_L;y)\\&+(\nu/2)(\sum\nolimits_{l=1}^{L}\Vert z_l-W_lp_l-b_l\Vert^2_2+\sum\nolimits_{l=1}^{L-1}\Vert q_l-f_l(z_l)\Vert^2_2)\\&s.t. \ p_{l+1}=q_l
\end{align*}
\end{problem}
where   $\textbf{p}=\{p_l\}_{l=1}^{L}$, $\textbf{W}=\{W_l\}_{l=1}^{L}$, $\textbf{b}=\{b_l\}_{l=1}^{L}$, $\textbf{z}=\{z_l\}_{l=1}^{L}$, $\textbf{q}=\{q_l\}_{l=1}^{L-1}$, and $\nu>0$ is a tuning parameter. As $\nu\rightarrow \infty$, Problem \ref{prob:problem 2} approaches Problem  \ref{prob:problem 1}. We reduce layer dependency by splitting the output of the $l$-th layer and the input of the $l+1$-th layer into two variables $p_{l+1}$ and $q_{l}$, respectively.

\begin{figure}
   \centering
    \includegraphics[width=\linewidth]{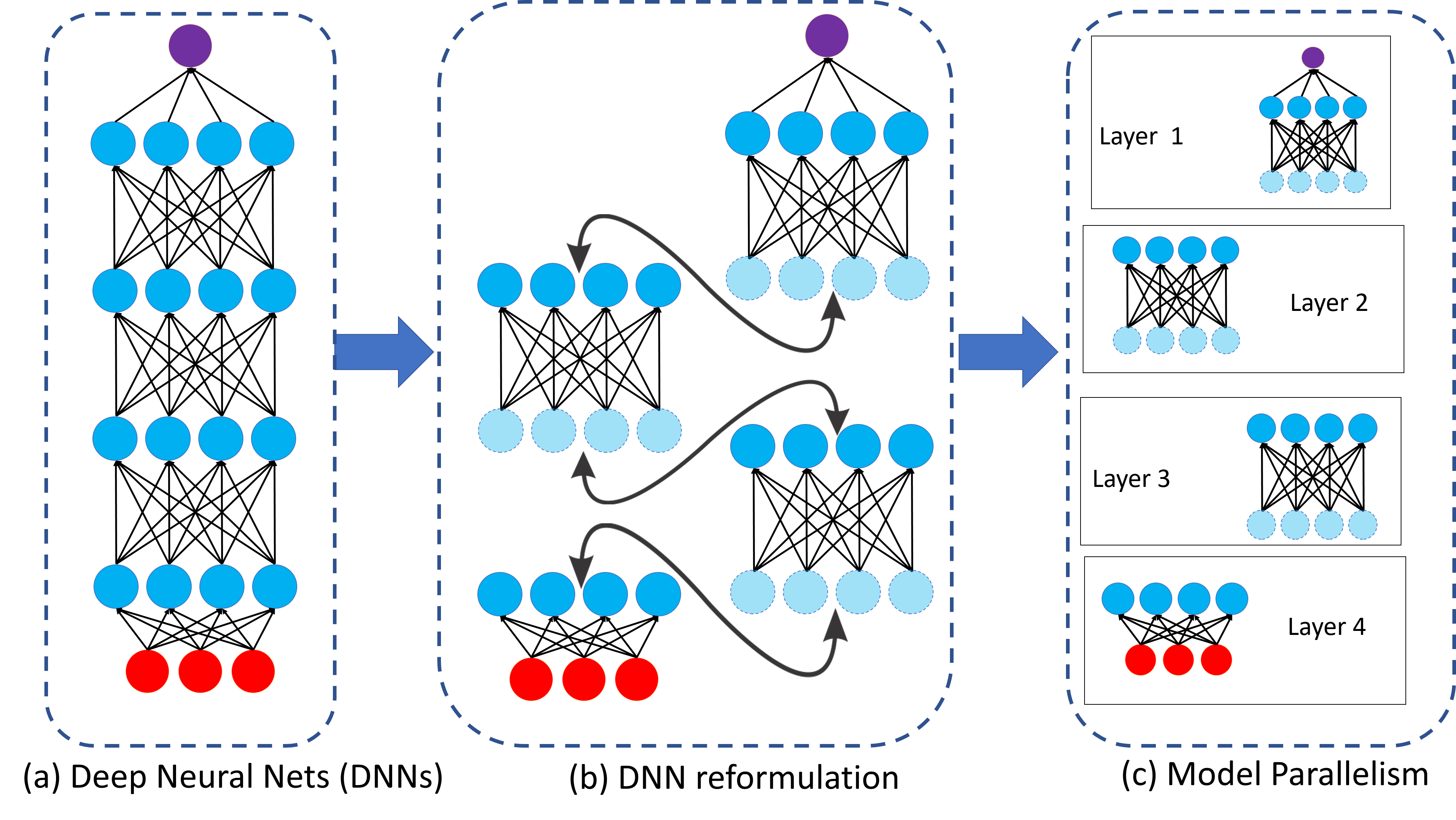}
    \caption{The pdADMM optimization framework: an overview}
    \label{fig:pdADMM framework}
\end{figure}

\indent  The high-level overview of the pdADMM algorithm is shown in Figure \ref{fig:pdADMM framework}.  Specifically, by breaking the whole neural network into multiple layers, each of which can be optimized by an independent worker. Therefore, the layerwise training can be implemented in parallel. Moreover, the gradient vanishing problem can be avoided in this way. This is because the accumulate gradient calculated by the backpropagation algorithm is split into layerwise components. \\
\indent Now we follow the ADMM routine to solve Problem \ref{prob:problem 2}, the augmented Lagrangian function is formulated mathematically as follows:
\begin{align*}
    &L_\rho(\textbf{p},\textbf{W},\textbf{b},\textbf{z},\textbf{q},\textbf{u})\\&=F(\textbf{p},\textbf{W},\textbf{b},\textbf{z},\textbf{q})+\sum\nolimits_{l\!=\!1}^{L\!-\!1}(u_l^T(p_{l\!+\!1}\!-\!q_l)\!+\!(\rho/2)\Vert p_{l+1}-q_l\Vert^2_2)\\&=\!R(z_L;y)\!+\!\phi(p_1,W_1,b_1,z_1)\!+\!\sum\nolimits_{l\!=\!2}^{L}\!\phi(p_l,W_l,b_l,z_l,q_{l-1},u_{l\!-\!1})\!\\&+\!(\nu/2)\sum\nolimits_{l\!=\!1}^{L-1}\Vert q_l\!-\!f_l(z_l)\Vert^2_2
\end{align*}
where $\phi(p_1,W_1,b_1,z_1)=(\nu/2)\Vert z_1-W_1p_1-b_1\Vert^2_2$, $\phi(p_l,W_l,b_l,z_l,q_{l-1},u_{l-1})=(\nu/2)\Vert z_l-W_lp_l-b_l\Vert^2_2+u^T_{l-1}(p_l-q_{l-1})+(\rho/2)\Vert p_l-q_{l-1}\Vert^2_2$,  $u_l(l=1,\cdots,L-1)$ are dual variables, $\rho>0$ is a parameter, and $\textbf{u}=\{u_l\}_{l=1}^{L-1}$.
The detail of the pdADMM is shown in Algorithm \ref{algo:distributed ADMM}. Specifically, Lines 5-9 update primal variables $\textbf{p}$, $\textbf{W}$, $\textbf{b}$, $\textbf{z}$ and $\textbf{q}$, respectively, while Line 11 updates the dual variable $\textbf{u}$. the discussion on how to solve subproblems generated by pdADMM is detailed in the next section.
\begin{algorithm} 
\caption{the pdADMM Algorithm} 
\begin{algorithmic}
\label{algo:distributed ADMM}
\REQUIRE $y$, $p_1=x$, $\rho$, $\nu$. 
\ENSURE $\textbf{p},\textbf{W},\textbf{b},\textbf{z},\textbf{q}$. 
\STATE Initialize $k=0$.
\WHILE{$\textbf{p}^k,\textbf{W}^{k},\textbf{b}^{k},\textbf{z}^{k},\textbf{q}^{k}$ not converged}
\STATE Update $p_l^{k+1}$ of different $l$ by Equation \eqref{eq:update p} in parallel.
\STATE Update $W_l^{k+1}$ of different $l$ by Equation \eqref{eq:update W} in parallel.
\STATE Update $b_l^{k+1}$ of different $l$ by Equation \eqref{eq:update b} in parallel.
\STATE Update $z_l^{k+1}$ of different $l$ by Equations \eqref{eq:update z} and \eqref{eq:update zl} in parallel.
\STATE Update $q_l^{k+1}$ of different $l$ by Equation \eqref{eq:update q} in parallel.
\STATE $r^k_l\leftarrow p^{k+1}_{l+1}-q^{k+1}_l(l=1,\cdots,L)$ in parallel \# Compute residuals.
\STATE Update $u_l^{k+1}$ of different $l$ by Equation \eqref{eq:update u} in parallel.
\STATE $k\leftarrow k+1$.
\ENDWHILE
\STATE Output $\textbf{p},\textbf{W},\textbf{b},\textbf{z},\textbf{q}$.
\end{algorithmic}
\end{algorithm}

\subsection{Solutions to All Subproblems}
\label{sec:subproblem}
\indent In this section, we discuss how to solve all subproblems generated by pdADMM in detail.\\
\textbf{1. Update $\textbf{p}^{k+1}$}\\
\indent The variable $\textbf{p}^{k+1}$ is updated as follows:
\begin{align*}
    p^{k+1}_l&\leftarrow \arg\min\nolimits_{p_l} L_\rho(\textbf{p},\textbf{W}^k,\textbf{b}^k,\textbf{z}^k,\textbf{q}^k,\textbf{u}^k)\\&=\phi(p_l,W^k_l,b^k_l,z^k_l,q^k_{l-1},u^k_{l-1})
\end{align*}
Because $W_l$ and $p_l$ are coupled in $\phi$, solving $p_l$ should require the time-consuming operation of matrix inversion of $W_l$. To handle this, we apply similar quadratic approximation techniques as used in dlADMM \cite{wang2019admm} as follows:
\begin{align}
p_l^{k+1}&\leftarrow\arg\min_{p_l} U_l(p_l;\tau_l^{k+1})
    \label{eq:update p}
\end{align}
where $
    U_l(p_l;\tau^{k+1}_l)=\phi(p^k_l,W^k_l,b^k_l,z^k_l,q^k_{l-1},u^k_{l-1})+(\nabla_{p^k_l}\phi(p^k_l,W^k_l,b^k_l,z^k_l,q^k_{l-1},u^k_{l-1})(p_l-p^k_l)+({\tau}^{k+1}_l/2)\Vert p_{l}-p^k_{l}\Vert^2_2
$, and $\tau^{k+1}_l>0$ is a parameter. $\tau_l^{k+1}$ should satisfy $\phi(p^{k+1}_l,W^k_l,b^k_l,z^k_l,q^k_{l-1},u^k_{l-1})\leq U_l(p^{k+1}_l;\tau^{k+1}_l)$.  The solution to Equation \eqref{eq:update p} is:
${p}^{k+1}_{l}\leftarrow p^k_{l}-\nabla_{p^k_{l}}\phi(p^k_l,W^k_l,b^k_l,z^k_l,q^k_{l-1},u^k_{l-1})/{\tau}^{k+1}_{l}$.\\
\textbf{2. Update $\textbf{W}^{k+1}$}\\
\indent The variable $\textbf{W}^{k+1}$ is updated as follows:
\begin{align*}
    W^{k+1}_l&\leftarrow \arg\min\nolimits_{W_l} L_\rho(\textbf{p}^{k+1},\textbf{W},\textbf{b}^k,\textbf{z}^{k},\textbf{q}^{k},\textbf{u}^{k})\\& =\arg\min\nolimits_{W_l} \begin{cases} \phi(p^{k+1}_1,W_1,b^k_1,z^k_1)
    &l=1\\\phi(p^{k\!+\!1}_l,W_l,b^k_l,z^k_l,q^k_{l\!-\!1},u^k_{l\!-\!1})&1<\!l\!\leq\! L
    \end{cases}
\end{align*}
Similar to updating $p_l$, the following subproblem should be solved instead:
\begin{align}
    W_l^{k+1}&\leftarrow\arg\min\nolimits_{W_l} V_l(W_l;\theta_l^{k+1}) \label{eq:update W}
\end{align}
where 
\begin{align*}
&V_1(W_1;\theta^{k+1}_1)= \phi(p^{k+1}_1,W^k_1,b^k_1,z^k_1)\\&+\!\nabla_{W^k_1}\phi^T(p^{k+1}_1,W^k_1,b^k_1,z^k_1)(W_1\!-\!W^k_1)   \\&+({\theta}^{k+1}_l/2)\Vert W_{1}-W^k_{1}\Vert^2_2\\&
 V_l(W_l;\theta^{k+1}_l)=
   \phi(p^{k+1}_l,W^k_l,b^k_l,z^k_l,q^k_{l-1},u^k_{l-1})\\&+\nabla_{W^k_l}\phi^T(p^{k+1}_l,W^k_l,b^k_l,z^k_l,q^k_{l-1},u^k_{l-1})(W_l-W^k_l)\\&+({\theta}^{k+1}_l/2)\Vert W_{l}-W^k_{l}\Vert^2_2
\end{align*}
and $\theta^{k+1}_l$ is a parameter, which should satisfy $\phi(p^{k+1}_1,W^{k+1}_1,b^k_1,z^k_1)\leq V(W^{k+1}_1;\theta^{k+1}_1)$ and $\phi(p^{k+1}_l,W^{k+1}_l,b^k_l,z^k_l,q^k_{l-1},u^k_{l-1})\leq V(W^{k+1}_l;\theta^{k+1}_l)(1<l<L)$. The solution to Equation \eqref{eq:update W} is shown as follows:
\begin{align*}
    {W}^{k\!+\!1}_{l}\!\leftarrow\! W^k_{l}\!-\!\begin{cases}\nabla_{W^k_1}\phi(p^{k+1}_1\!,\!W^k_1\!,\!b^k_1\!,\!z^k_1)/\theta^{k+1}_l& l=1\\\nabla_{W^k_l}\phi(p^{k\!+\!1}_l\!,\!W^k_l\!,\!b^k_l\!,\!z^k_l\!,\!q^k_{l\!-\!1}\!,\!u^k_{l\!-\!1})/\theta^{k\!+\!1}_l& 1\!<l\!\leq\! L
    \end{cases}
\end{align*}
\textbf{3. Update $\textbf{b}^{k+1}$}\\
\indent The variable $\textbf{b}^{k+1}$ is updated as follows:
\begin{align*}
    b^{k+1}_l&\leftarrow \arg\min\nolimits_{b_l} L_\rho(\textbf{p}^{k+1},\textbf{W}^{k+1},\textbf{b},\textbf{z}^{k},\textbf{q}^{k},\textbf{u}^{k})\\& =\arg\min\nolimits_{b_l} \begin{cases} \phi(p^{k+1}_1,W^{k+1}_1,b_1,z^k_1)
    &l\!=\!1\\\phi(p^{k+1}_l,W^{k\!+\!1}_l,b_l,z^k_l,q^k_{l\!-\!1},u^k_{l\!-\!1})&1<\!l\!\leq L
    \end{cases}.
\end{align*}
Similarly, we solve the following subproblems instead:
\begin{align*}
&b^{k+1}_1\leftarrow \arg\min\nolimits_{b_1} \phi(p^{k+1}_1,W^{k+1}_1,b^k_1,z^k_1)\\&+\nabla_{b^k_1}\phi^T(p^{k+1}_1,W^{k+1}_1,b^k_1,z^k_1)(b_l-b^k_l)+(\nu/2)\Vert b_{l}-b^k_{l}\Vert^2_2
\end{align*}
\begin{align}
&\nonumber b^{k+1}_l\leftarrow \arg\min\nolimits_{b_l} \phi(p^{k+1}_l,W^{k+1}_l,b^k_l,z^k_l,q^k_{l-1},u^k_{l-1})\\&\nonumber+\nabla_{b^k_l}\phi^T(p^{k+1}_l,W^{k+1}_l,b^k_l,z^k_l,q^k_{l-1},u^k_{l-1})(b_l-b^k_l)\\&+(\nu/2)\Vert b_{l}-b^k_{l}\Vert^2_2 (1< l\leq L) \label{eq:update b}
\end{align}
\indent The solution to Equation \eqref{eq:update b} is:
\begin{align*}
    {b}^{k\!+\!1}_{l}\leftarrow b^k_{l}\!-\!\begin{cases}\nabla_{b^k_1}\phi(p^{k\!+\!1}_1,W^{k\!+\!1}_1,b^k_1,z^k_1)/\nu& l=1\\\nabla_{b^k_l}\phi(p^{k\!+\!1}_l,W^{k\!+\!1}_l,b^k_l,z^k_l,q^k_{l\!-\!1},u^k_{l\!-\!1})/\nu& 1\!<\!l\!\leq L
    \end{cases}
\end{align*}
\textbf{4. Update $\textbf{z}^{k+1}$}\\
\indent The variable
$\textbf{z}^{k+1}$ is updated as follows:
\begin{align}
    \nonumber& z^{k\!+\!1}_l\!\leftarrow\! \arg\min\nolimits_{z_l} (\nu/2)\Vert z_l\!-\!W^{k\!+\!1}_l\!p^{k\!+\!1}_l\!-\!b^{k\!+\!1}_l\Vert^2_2\\&+\!(\nu/2)\Vert q^k_l\!-
    \!f_l(z_l)\Vert^2_2\!+\!(\nu/2)\Vert z_l\!-\!z^k_l\Vert^2_2(l\!<\! L)\label{eq:update z}\\
    &z^{k\!+\!1}_L\!\leftarrow\! \arg\min\nolimits_{z_l}\! R(z_L;y)\!+\!(\nu/2)\!\Vert z_L\!-\!W^{k\!+\!1}_Lp^{k\!+\!1}_L\!-\!b^{k\!+\!1}_L\Vert^2_2\! \label{eq:update zl}
\end{align}
where a quadratic term $(\nu/2)\Vert z_l-z^k_l\Vert^2_2$ is added in Equation \eqref{eq:update z} to control $z^{k+1}_l$ to close to $z^{k}_l$.
Equation \eqref{eq:update zl} is convex, which can be solved by Fast Iterative Soft Thresholding Algorithm (FISTA) \cite{beck2009fast}.\\
For Equation \eqref{eq:update z}, nonsmooth activations usually lead to  closed-form solutions \cite{wang2019admm,wang2017nonconvex}. For example, for Relu $f_l(z_l)=\max(z_l,0)$, the solution to Equation \eqref{eq:update z} is shown as follows:
\begin{align*}
    z^{k+1}_l=
    \begin{cases} \min((W^{k\!+\!1}_{l}p^{k\!+\!1}_{l}\!+\!b^{k+1}_l\!+\!z^k_l)/2,0)&z^{k\!+\!1}_l\leq0\\
    \max((W^{k\!+\!1}_{l}p^{k\!+\!1}_{l}\!+\!b^{k\!+\!1}_l\!+\!q^k_l\!+\!z^k_l)/3,0)&z^{k+1}_l\geq0\\
    \end{cases}
\end{align*}
For smooth activations such as tanh and sigmoid, a lookup-table is recommended \cite{wang2019admm}.\\
\textbf{5. Update $\textbf{q}^{k+1}$}\\
\indent The variable
$\textbf{q}^{k+1}$ is updated as follows:
\begin{align}
     \nonumber q^{k+1}_l&\leftarrow \arg\min\nolimits_{q_l} L_\rho(\textbf{p}^{k+1},\textbf{W}^{k+1},\textbf{b}^{k+1},\textbf{z}^{k+1},\textbf{q},\textbf{u}^k)\\&=\arg\min\nolimits_{q_l}\phi(p^{k+1}_{l+1},W^{k+1}_{l+1},b^{k+1}_{l+1},z^{k+1}_{l+1},q_l,u^{k}_l). \label{eq:update q}
\end{align}
Equation \eqref{eq:update q} has a closed-form solution as follows:\\
\begin{align*}
     q^{k+1}_l\leftarrow (\rho p_{l+1}^{k+1}+u^k_l+\nu f_l(z^{k+1}_l))/(\rho+\nu)
\end{align*}
\textbf{6. Update $\textbf{u}^{k+1}$}\\
\indent The variable
$\textbf{u}^{k+1}$ is updated as follows:
\begin{align}
u^{k+1}_l&\leftarrow u^k_l+\rho(p^{k+1}_{l+1}-q^{k+1}_{l}) \label{eq:update u}
\end{align}
\indent Finally, Our proposed pdADMM can be efficient for training a deep feed-forward neural network. To achieve this,  we begin from training a swallow neural network with the first few layers of the deep neural network, then more layers are added for training step by step until finally all layers are involved in the training process. The pdADMM can achieve  good performance as well as reduce training cost by this strategy.
\section{Convergence Analysis}
\label{sec:convergence}
\indent In this section, the theoretical convergence of the proposed pdADMM algorithm. Firstly, the Lipschitz continuity and coercivity are defined as follows:
\begin{definition}(Lipschitz Continuity)
A function $g(x)$ is Lipschitz continuous if there exists a constant $D>0$ such that $\forall x_1,x_2$, the following holds
\begin{align*}
    \Vert g(x_1)-g(x_2)\Vert\leq D\Vert x_1-x_2\Vert. 
\end{align*}
\end{definition}
\begin{definition}(Coercivity)
A function $h(x)$ is coerce over the feasible set  $\mathscr{F}$ means that $h(x)\rightarrow \infty$ if $x\in \mathscr{F}$ and $\Vert x\Vert\rightarrow \infty$.
\end{definition}
Then the following assumption is required for convergence analysis.
\begin{assumption}
$f_l(z_l)$ is Lipschitz continuous with coefficient $S>0$, and $F(\textbf{p},\textbf{W},\textbf{b},\textbf{z},\textbf{q})$ is coercive. Moreover, $\partial f_l(z_l)$ is bounded, i.e. there exists $M>0$ such that $\Vert\partial f_l(z_l)\Vert \leq M$. 
\label{ass:lipschitz continuous}
\end{assumption}
Assumption \ref{ass:lipschitz continuous} is mild to satisfy: most common activation functions such as Relu and leaky Relu satisfy Assumption \ref{ass:lipschitz continuous}.
No assumption is needed on the risk function $R(z_l;y)$, which shows that the convergence condition of our proposed pdADMM is milder than that of the dlADMM, which requires $R(z_l;y)$ to be Lipschitz differentiable \cite{wang2019admm}. Due to space limit, the detailed proofs are provided in the Appendix\footnote{Proofs: \url{https://github.com/xianggebenben/Junxiang_Wang/blob/master/supplementary_material/ICDM2020/pdADMM.pdf}.}. The technical proofs follow the similar routine as dlADMM \cite{wang2019admm}. The difference consists in the fact that the dual variable $u_l$ is controlled by $q_l$ and $z_l$ (Lemma \ref{lemma:u square bound} in the Appendix), which holds under Assumption \ref{ass:lipschitz continuous}, while $u_l$ can  be controlled only by $z_l$ in the convergence proof of dlADMM. The first lemma shows that the objective keeps decreasing when $\rho$ is sufficently large.

\begin{lemma}[Decreasing Objective]
\label{lemma:objective decrease}
If $\rho>\max(4\nu S^2,(\sqrt{17}+1)\nu/2)$, there exist $C_1=\nu/2-2\nu^2S^2/\rho>0$ and $C_2=\rho/2-2\nu^2/\rho-\nu/2>0$ such that it holds for any $k\in \mathbb{N}$ that
\begin{align}
    \nonumber &L_\rho(\textbf{p}^{k},\textbf{W}^{k},\textbf{b}^k,\textbf{z}^{k},\textbf{q}^k,\textbf{u}^{k})\!-\!L_\rho(\textbf{p}^{k\!+\!1},\textbf{W}^{k\!+\!1},\textbf{b}^{k\!+\!1},\textbf{z}^{k\!+\!1},\textbf{q}^{k\!+\!1},\textbf{u}^{k\!+\!1})\\\nonumber &\geq \sum\nolimits_{l\!=\!2}^L (\tau^{k\!+\!1}_l/2)\Vert p^{k+1}_l\!-\!p^k_l\Vert^2_2\!+\!\sum\nolimits_{l\!=\!1}^{L}(\theta^{k\!+\!1}_l/2)\Vert  W^{k\!+\!1}_l\!-\!W^k_l\Vert^2_2\!\\&\nonumber+\!\sum\nolimits_{l\!=\!1}^{L}(\nu/2)\Vert  b^{k\!+\!1}_l\!-\!b^k_l\Vert^2_2\!+\!\sum\nolimits_{l\!=\!1}^{L\!-\!1} C_1\Vert z^{k+1}_l\!-\!z^k_l\Vert^2_2\\&\!+\!(\nu/2)\Vert z^{k+1}_L-z^k_L\Vert^2_2+\sum\nolimits_{l=1}^{L-1}C_2\Vert q^{k+1}_l-q^k_l\Vert^2_2
    \label{eq:objective decrease}
\end{align}
\end{lemma}
\indent The second Lemma illustrates that the objective is bounded from below when $\rho$ is large enough, and all variables are bounded.
\begin{lemma} [Bounded Objective]
\label{lemma:lower bounded}
If $\rho> \nu$, then
 $L_\rho(\textbf{p}^k,\textbf{W}^{k},\textbf{b}^{k},\textbf{z}^{k},\textbf{q}^{k},\textbf{u}^{k})$
 is lower bounded. Moreover, $\textbf{p}^k,\textbf{W}^k,\textbf{b}^k,\textbf{z}^k,\textbf{q}^{k}$,and $\textbf{u}^k$ are bounded, i.e. there exist $\mathbb{N}_\textbf{p}$, $\mathbb{N}_\textbf{W}$, $\mathbb{N}_\textbf{b}$,  $\mathbb{N}_\textbf{z}$,  $\mathbb{N}_\textbf{q}$, and $\mathbb{N}_\textbf{u}>0$, such that $\Vert \textbf{p}^k\Vert\leq \mathbb{N}_\textbf{p}$, $\Vert \textbf{W}^k\Vert\leq \mathbb{N}_\textbf{W}$, $\Vert \textbf{b}^k\Vert\leq \mathbb{N}_\textbf{b}$, $\Vert \textbf{z}^k\Vert\leq \mathbb{N}_\textbf{z}$, $\Vert \textbf{q}^k\Vert\leq \mathbb{N}_\textbf{q}$, and $\Vert \textbf{u}^k\Vert\leq \mathbb{N}_\textbf{u}$. \end{lemma}
 \indent Based on Lemmas \ref{lemma:objective decrease} and \ref{lemma:lower bounded}, the following theorem ensures that the objective  is convergent.
 \begin{theorem}[Convergent Objective] \label{theo:convergent variable}
If $\rho>\max(4\nu S^2,(\sqrt{17}+1)\nu/2)$, then
$L_\rho(\textbf{p}^k,\textbf{W}^k,\textbf{b}^k,\textbf{z}^k,\textbf{q}^{k},\textbf{u}^k)$ is convergent. Moreover, $\lim_{k\rightarrow\infty}\Vert\textbf{p}^{k+1}-\textbf{p}^{k}\Vert^2_2=0$, $\lim_{k\rightarrow\infty}\Vert\textbf{W}^{k+1}-\textbf{W}^{k}\Vert^2_2=0$, $\lim_{k\rightarrow\infty}\Vert\textbf{b}^{k+1}-\textbf{b}^{k}\Vert^2_2=0$, $\lim_{k\rightarrow\infty}\Vert\textbf{z}^{k+1}-\textbf{z}^{k}\Vert^2_2=0$, $\lim_{k\rightarrow\infty}\Vert\textbf{q}^{k+1}-\textbf{q}^{k}\Vert^2_2=0$, $\lim_{k\rightarrow\infty}\Vert\textbf{u}^{k+1}-\textbf{u}^{k}\Vert^2_2=0$.
\end{theorem}
\begin{proof}
From Lemmas \ref{lemma:objective decrease} and \ref{lemma:lower bounded}, we know that $L_\rho(\textbf{p}^{k},\textbf{W}^{k},\textbf{b}^k,\textbf{z}^{k},\textbf{q}^k,\textbf{u}^{k})$ is convergent because a monotone bounded sequence converges. Moreover, we take the limit on the both sides of Inequality \eqref{eq:objective decrease} to obtain
\begin{align*}
     &0=\lim\nolimits_{k\rightarrow\infty}L_\rho(\textbf{p}^{k},\textbf{W}^{k},\textbf{b}^k,\textbf{z}^{k},\textbf{q}^k,\textbf{u}^{k})\\&-\lim\nolimits_{k\rightarrow\infty}L_\rho(\textbf{p}^{k+1},\textbf{W}^{k+1},\textbf{b}^{k+1},\textbf{z}^{k+1},\textbf{q}^{k+1},\textbf{u}^{k+1})\\ &\geq \lim\nolimits_{k\rightarrow\infty}(\sum\nolimits_{l=2}^L (\tau^{k+1}_l/2)\Vert p^{k+1}_l\!-\!p^k_l\Vert^2_2\\&+\sum\nolimits_{l\!=\!1}^{L}(\theta^{k\!+\!1}_l/2)\Vert  W^{k\!+\!1}_l\!-\!W^k_l\Vert^2_2+\sum\nolimits_{l\!=\!1}^{L}(\nu/2)\Vert  b^{k\!+\!1}_l\!-\!b^k_l\Vert^2_2\\&+\sum\nolimits_{l=1}^{L-1} C_1\Vert z^{k+1}_l-z^k_l\Vert^2_2+(\nu/2)\Vert z^{k+1}_L-z^k_L\Vert^2_2\\&+\sum\nolimits_{l=1}^{L-1}C_2\Vert q^{k+1}_l-q^k_l\Vert^2_2)\geq 0
\end{align*}
Because $L_\rho(\textbf{p}^{k},\textbf{W}^{k},\textbf{b}^k,\textbf{z}^{k},\textbf{q}^k,\textbf{u}^{k})$ is convergent, then $\lim_{k\rightarrow\infty}\Vert\textbf{p}^{k+1}-\textbf{p}^{k}\Vert^2_2=0$, $\lim_{k\rightarrow\infty}\Vert\textbf{W}^{k+1}-\textbf{W}^{k}\Vert^2_2=0$, $\lim_{k\rightarrow\infty}\Vert\textbf{b}^{k+1}-\textbf{b}^{k}\Vert^2_2=0$, $\lim_{k\rightarrow\infty}\Vert\textbf{z}^{k+1}-\textbf{z}^{k}\Vert^2_2=0$, and $\lim_{k\rightarrow\infty}\Vert\textbf{q}^{k+1}-\textbf{q}^{k}\Vert^2_2=0$. $\lim_{k\rightarrow\infty}\Vert\textbf{u}^{k+1}-\textbf{u}^{k}\Vert^2_2=0$ is derived from Lemma \ref{lemma:u square bound} in the Appendix.
\end{proof}
 \indent The third lemma guarantees that the subgradient of the objective is upper bounded, which is stated as follows:
 \begin{lemma}[Bounded Subgradient]
\label{lemma:subgradient bound}
There exists a constant $C>0$ and $g^{k+1}\in \partial L_\rho(\textbf{p}^{k+1},\textbf{W}^{k+1},\textbf{b}^{k+1},\textbf{z}^{k+1},\textbf{q}^{k+1},\textbf{u}^{k+1})$ such that
\begin{align*}
    &\Vert g^{k+1}\Vert \leq C(\Vert\textbf{p}^{k+1}-\textbf{p}^{k}\Vert+\Vert\textbf{W}^{k+1}-\textbf{W}^{k}\Vert+\Vert\textbf{b}^{k+1}-\textbf{b}^{k}\Vert\\&+\Vert\textbf{z}^{k+1}-\textbf{z}^{k}\Vert+\Vert\textbf{q}^{k+1}-\textbf{q}^{k}\Vert+\Vert\textbf{u}^{k+1}-\textbf{u}^{k}\Vert)
\end{align*}
\end{lemma}
 \indent Now based on Theorem \ref{theo:convergent variable}, and Lemma \ref{lemma:subgradient bound}, the convergence of the pdADMM algorithm to a critical point is presented in the following theorem.
\begin{theorem} [Convergence to a Critical Point]
\label{theo: global convergence}
If $\rho>\max(4\nu S^2,(\sqrt{17}+1)\nu/2)$, then for the variables $(\textbf{p},\textbf{W},\textbf{b},\textbf{z},\textbf{q},\textbf{u})$ in Problem \ref{prob:problem 2}, starting from any $(\textbf{p}^{0},\textbf{W}^{0},\textbf{b}^{0},\textbf{z}^{0},\textbf{q}^{0},\textbf{u}^{0})$, $(\textbf{p}^{k},\textbf{W}^{k},\textbf{b}^{k},\textbf{z}^{k},\textbf{q}^{k},\textbf{u}^{k})$ has at least a limit point $(\textbf{p}^*,\textbf{W}^*,\textbf{b}^*,\textbf{z}^*,\textbf{q}^*,\textbf{u}^*)$, and any limit point is a critical point of Problem \ref{prob:problem 2}. That is, $0\in \partial L_\rho(\textbf{p}^*,\textbf{W}^*,\textbf{b}^*,\textbf{z}^*,\textbf{q}^*,\textbf{u}^*)$. In other words, 
\begin{align*}
    & p^*_{l+1}=q^*_l, \ \nabla_{\textbf{p}^*} L_\rho(\textbf{p}^*,\textbf{W}^*,\textbf{b}^*,\textbf{z}^*,\textbf{q}^*,\textbf{u}^*)=0,\\& \
    \nabla_{\textbf{W}^*} L_\rho(\textbf{p}^*,\textbf{W}^*,\textbf{b}^*,\textbf{z}^*,\textbf{q}^*,\textbf{u}^*)=0,\\& \nabla_{\textbf{b}^*} L_\rho(\textbf{p}^*,\textbf{W}^*,\textbf{b}^*,\textbf{z}^*,\textbf{q}^*,\textbf{u}^*)=0, \\& 0\in\partial_{\textbf{z}^*} L_\rho(\textbf{p}^*,\textbf{W}^*,\textbf{b}^*,\textbf{z}^*,\textbf{q}^*,\textbf{u}^*),\\& 
     \nabla_{\textbf{q}^*} L_\rho(\textbf{p}^*,\textbf{W}^*,\textbf{b}^*,\textbf{z}^*,\textbf{q}^*,\textbf{u}^*)=0.
\end{align*}
\end{theorem}
\begin{proof}
 From Lemma \ref{lemma:lower bounded}, $(\textbf{p}^k,\textbf{W}^k, \textbf{b}^k, \textbf{z}^k, \textbf{q}^k,\textbf{u}^k)$ has at least a limit point $(\textbf{p}^*,\textbf{W}^*, \textbf{b}^*, \textbf{z}^*, \textbf{q}^*,\textbf{u}^*)$ because a bounded sequence has at least a limit point. From Lemma \ref{lemma:subgradient bound} and Theorem \ref{theo:convergent variable},
$\Vert g^{k+1}\Vert \rightarrow 0$ as $k\rightarrow \infty$. According to the definition of general subgradient (Defintion 8.3 in \cite{rockafellar2009variational}), we have $0\in \partial L_\rho(\textbf{p}^*,\textbf{W}^*,\textbf{b}^*,\textbf{z}^*,\textbf{q}^*,\textbf{u}^*)$. In other words, every limit point $(\textbf{p}^*,\textbf{W}^*, \textbf{b}^*, \textbf{z}^*, \textbf{q}^*,\textbf{u}^*)$ is a critical point.
\end{proof}
\indent Theorem \ref{theo: global convergence} shows that our proposed pdADMM algorithm converges for sufficiently large $\rho$, which is consistent with previous literature \cite{wang2019admm}. Next, the following theorem ensures the sublinear convergence rate $o(1/k)$ of the proposed pdADMM algorithm, whose proof is at the end of this paper.
\begin{theorem}[Convergence Rate]
For a sequence $(\textbf{p}^k,\textbf{W}^k,\textbf{b}^k,\textbf{z}^k,\textbf{q}^k,\textbf{u}^k)$, define $c_k=\min\nolimits_{0\leq i\leq k}(\sum\nolimits_{l=2}^L (\tau^{i+1}_l/2)\Vert p^{i+1}_l\!-\!p^i_l\Vert^2_2+\sum\nolimits_{l\!=\!1}^{L}(\theta^{i\!+\!1}_l/2)\Vert  W^{i\!+\!1}_l\!-\!W^i_l\Vert^2_2+\sum\nolimits_{l\!=\!1}^{L}(\nu/2)\Vert  b^{i\!+\!1}_l\!-\!b^i_l\Vert^2_2+\sum\nolimits_{l=1}^{L-1} C_1\Vert z^{i+1}_l-z^i_l\Vert^2_2+(\nu/2)\Vert z^{i+1}_L-z^i_L\Vert^2_2+\sum\nolimits_{l=1}^{L-1}C_2\Vert q^{i+1}_l-q^i_l\Vert^2_2)$ where $C_1=\nu/2-2\nu^2S^2/\rho>0$ and $C_2=\rho/2-2\nu^2/\rho-\nu/2>0$, then the convergence rate of $c_k$ is $o(1/k)$.
\label{theo: convergence rate}
\end{theorem}
\section{Experiments}
\label{sec:experiment}
In this section, we evaluate the performance of the proposed pdADMM using six benchmark datasets. Speedup, convergence and accuracy performance are compared with several state-of-the-art optimizers.  All experiments were conducted on 64-bit machine  with Intel Xeon(R) silver 4114 Processor and 48GB RAM.

\subsection{Datasets}
 In this experiment, six benchmark datasets were used for performance evaluation:\\
 1. MNIST \cite{lecun1998gradient}. The MNIST dataset has ten classes of handwritten-digit images, which was firstly introduced by Lecun et al. in 1998 \cite{lecun1998gradient}. It contains 55,000 training samples and 10,000 test samples with 196 features each, which is provided by the Keras library \cite{chollet2015keras}. \\
2. Fashion MNIST \cite{xiao2017fashion}. The Fashion MNIST dataset has ten classes of assortment images on the website of Zalando, which is Europe’s largest online fashion platform \cite{xiao2017fashion}. The Fashion-MNIST dataset consists of 60,000 training samples and 10,000 test samples with 196 features each.\\
3. kMNIST(Kuzushiji-MNIST) \cite{clanuwat2018deep}.  The kMNIST dataset has ten classes, each of which is a character to represent each of the 10 rows of Hiragana. The kMNIST dataset consists of 60,000 training samples and 10,000 test samples with 196 features each.\\
4. SVHN (Street View House Numbers) \cite{SVHN}. The SVHN dataset is obtained from house numbers in Google Street View images. It consists of ten classes of digits. In our experiments, we use three classes '0', '1' and '2'. The number of training data and test data are 24,446 and 9,248, respectively, with 768 features each.\\
5. CIFAR10 \cite{krizhevsky2009learning}. CIFAR10 is a collection of color images with 10 different classes. In our experiments, we use two classes '0' and '6'. The number of training data and test data are 12,000 and 2,000, respectively, with 768 features each.\\
6. CIFAR100 \cite{krizhevsky2009learning}. CIFAR100 is similar to CIFAR10 except that CIFAR100 has 100 classes. In our experiments, we use two classes '0' and '2'. The number of training data and test data are 5,000 and 1,000, respectively, with 768 features each.
 \subsection{Speedup}
    \begin{figure}
    \centering
       \includegraphics[width=0.8\linewidth]{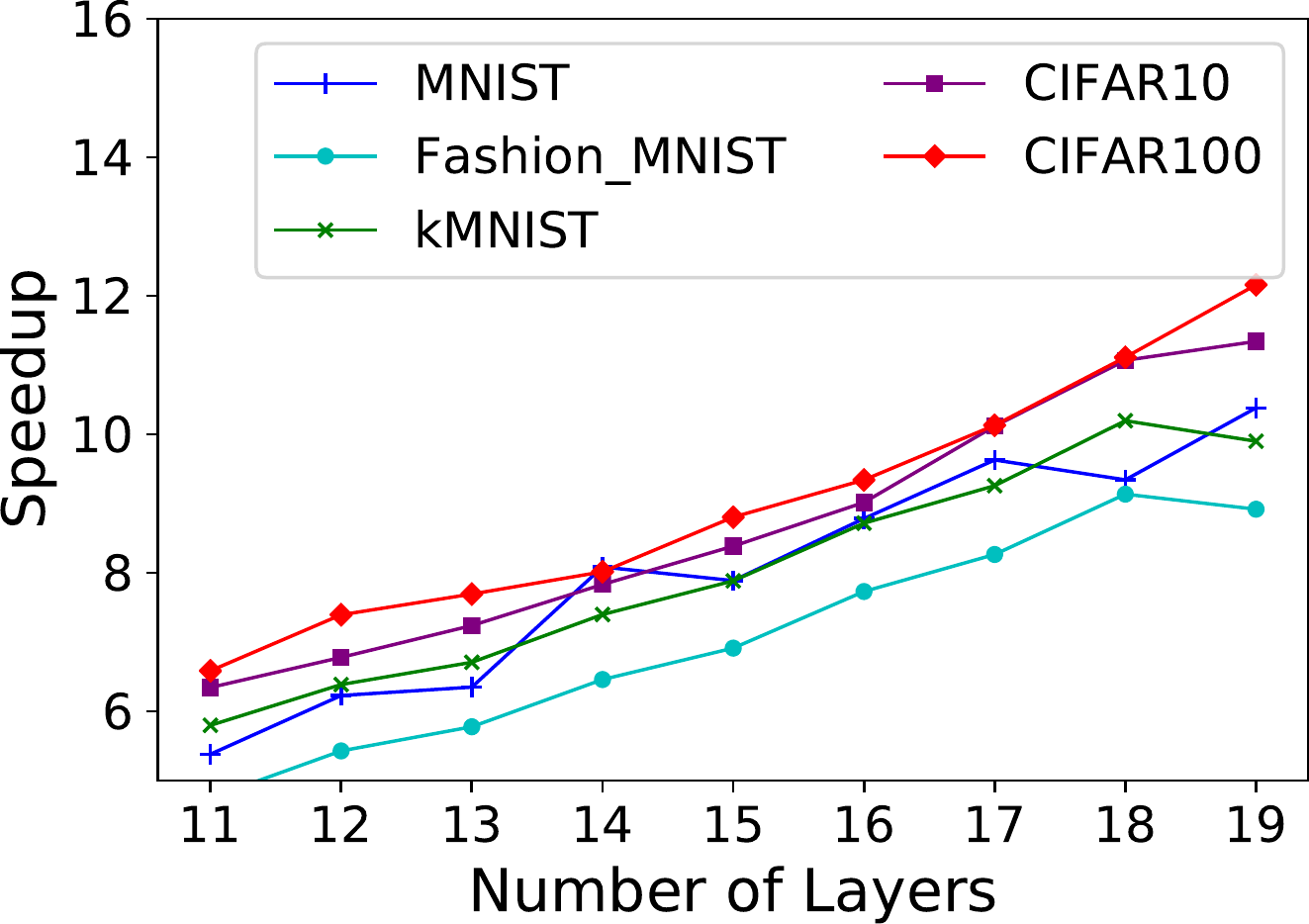}
       \caption{The relationship between speedup and the number of layers: the speedup increases linearly with the number of layers.}
       \label{fig:speedup}
       \vspace{-0.5cm}
\end{figure}

\begin{table}[]
\centering
\begin{tabular}{c|c|c|c}
\hline\hline
\multicolumn{4}{c}{MNIST dataset}\\\hline
 Neurons$\#$ & Serial pdADMM (sec)& pdADMM(sec)& Speedup\\\hline
1500&    237.66&    26.78&    8.87\\\hline
1600&    348.70&    31.78&    10.97\\\hline
1700&    390.51&    35.79&    10.91\\\hline
1800&    475.60&    41.37&    11.50\\\hline
1900&    465.57&    45.87&    10.15\\\hline
2000&    570.90&    50.70&    11.26\\
\hline
2000&    570.9&    50.7&    11.26\\\hline
2100&    570&    54.91&    10.38\\\hline
2200&    678.83&    63.59&    10.68\\\hline
2300&    710.3&    70.36&    10.10\\\hline
2400&    766.82&    62.5&    12.27\\
\hline\hline
\multicolumn{4}{c}{Fashion MNIST dataset}\\\hline
 Neurons$\#$ & Serial pdADMM (sec)& pdADMM(sec)& Speedup\\\hline
1500&    358.68&    32.65&    10.99\\\hline
1600&    407.71&    37.90&    10.76\\\hline
1700&    476.79&    44.75&    10.65\\\hline
1800&    539.51&    50.50&    10.68\\\hline
1900&    599.42&    53.88&    11.13\\\hline
2000&    645.87&    58.68&    11.01\\\hline
2100&    740.39&    67.91&    10.90\\\hline
2200&    818.58&    74.17&    11.03\\\hline\hline
\multicolumn{4}{c}{kMNIST dataset}\\\hline
 Neurons$\#$ & Serial pdADMM (sec)& pdADMM(sec)& Speedup\\\hline
 1500&    354.85&    32.65&    10.87\\\hline
1600&    407.73&    37.11&    10.99\\\hline
1700&    472.4648&    42.58&    11.10\\\hline
1800&    539.52&    48.78&    11.06\\\hline
1900&    596.84&    55.56&    10.74\\\hline
2000&    660.58&    56.10&    11.78\\\hline
2100&    737.78&    66.95&    11.02\\\hline
2200&    806.74&    76.16&    10.59\\
\hline\hline
\multicolumn{4}{c}{CIFAR10  dataset}\\\hline
 Neurons$\#$ & Serial pdADMM (sec)& pdADMM(sec)& Speedup\\\hline
 1500&    326.62&    25.00&    13.06\\\hline
1600&    374.82&    28.96&    12.94\\\hline
1700&    433.46&    33.99&    12.75\\\hline
1800&    485.86&    38.66&    12.57\\\hline
1900&    544.11&    43.10&    12.62\\\hline
2000&    572.33&    46.90&    12.20\\\hline
2100&    602.65&    55.25&    10.91\\\hline
2200&    732.79&    59.27&    12.36\\\hline
2300&    784.87&    56.26&    13.95\\\hline
2400&    854.47&    63.1&    13.54\\\hline
\hline
\multicolumn{4}{c}{CIFAR100  dataset}\\\hline
 Neurons$\#$ & Serial pdADMM (sec)& pdADMM(sec)& Speedup\\\hline
 1500&    334.55&    25.39&    13.18\\\hline
1600&    382.24&    29.3&    13.05\\\hline
1700&    445.23&    34&    13.09\\\hline
1800&    500.00&    38.38&    13.03\\\hline
1900&    549.77&    43.25&    12.71\\\hline
2000&    576.10&    42.47&    13.56\\\hline
2100&    666.06&    47.43&    14.04\\\hline
2200&    735.63&    52.41&    14.04\\\hline
2300&    793.03&    56.73&    13.98\\\hline
2400&    857.41&    62.3&    13.76\\\hline
\hline
\end{tabular}
\caption{The relation between speedup and number of  neurons on the MNIST, Fashion MNIST datasets,  kMNIST, CIFAR10 and CIFAR100 datasets: the pdADMM runs 10 times faster than its serial version.}
\label{tab: running time}
\end{table}
 In this experiment, we investigate the speedup of the proposed pdADMM algorithm concerning  the number of layers and the number of neurons on the large-scale deep neural networks. The activation function was set to the Rectified linear unit (Relu). The loss function was the cross-entropy loss. The running time per epoch was the average of 10 epochs. $\rho$ and $\nu$ were both set to $10^{-4}$.\\
  \indent Firstly, we investigated the relationship between speedup and the number of layers. We set up a feed-forward neural network with different number of hidden layers, which ranges from 11 to 19. The number of neurons in each layer was fixed to 2,400. The SVHN dataset was not tried due to memory issues. Figure \ref{fig:speedup} shows that the speedup increases linearly with the number of layers. Specifically, the speedup reached 11 when 19 hidden layers were trained.  This indicates that the deeper a neural network is, the more speedup our proposed pdADMM can gain.\\
 \indent Secondly, the relationship between speedup and the number of neurons was studied. Specifically, we test our proposed pdADMM algorithm on a feed-forward neural network with 19 hidden layers. The number of neurons in each layer ranges from 1,500 to 2,400. The speedup was shown in Table \ref{tab: running time} on the MNIST and Fashion MNIST datasets. Specifically, the speedup remains stable around $10$ no matter how many neurons were trained. This concludes that the speedup of the proposed pdADMM is independent of the number of neurons.
 \subsection{Convergence}
 \begin{figure*}
 \begin{minipage}{0.49\linewidth}
      \centering
     \includegraphics[width=0.6\linewidth]{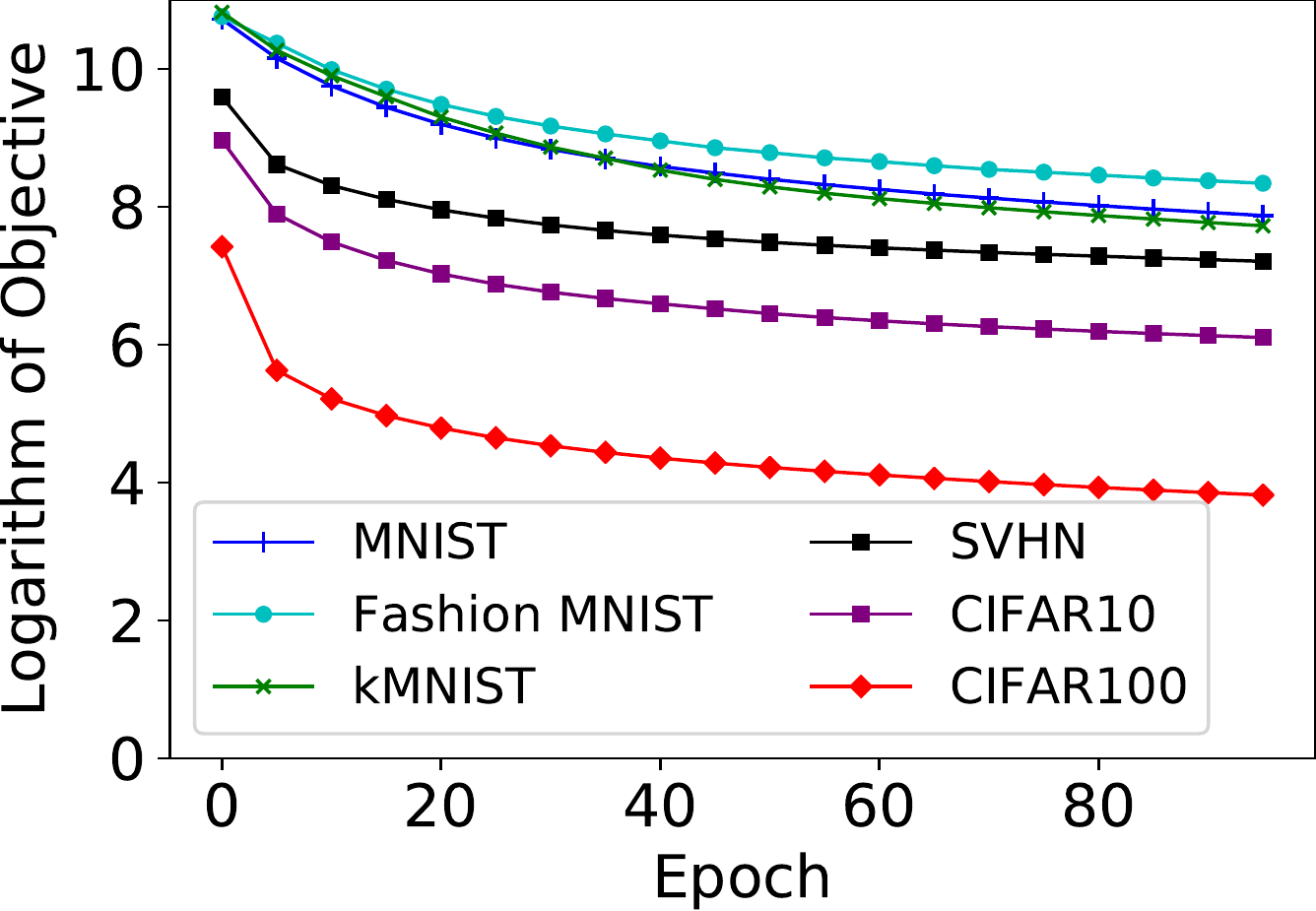}
     \centerline{(a). Objective versus epoch}
      \end{minipage}
      \hfill
 \begin{minipage}{0.49\linewidth}
      \centering
     \includegraphics[width=0.6\linewidth]{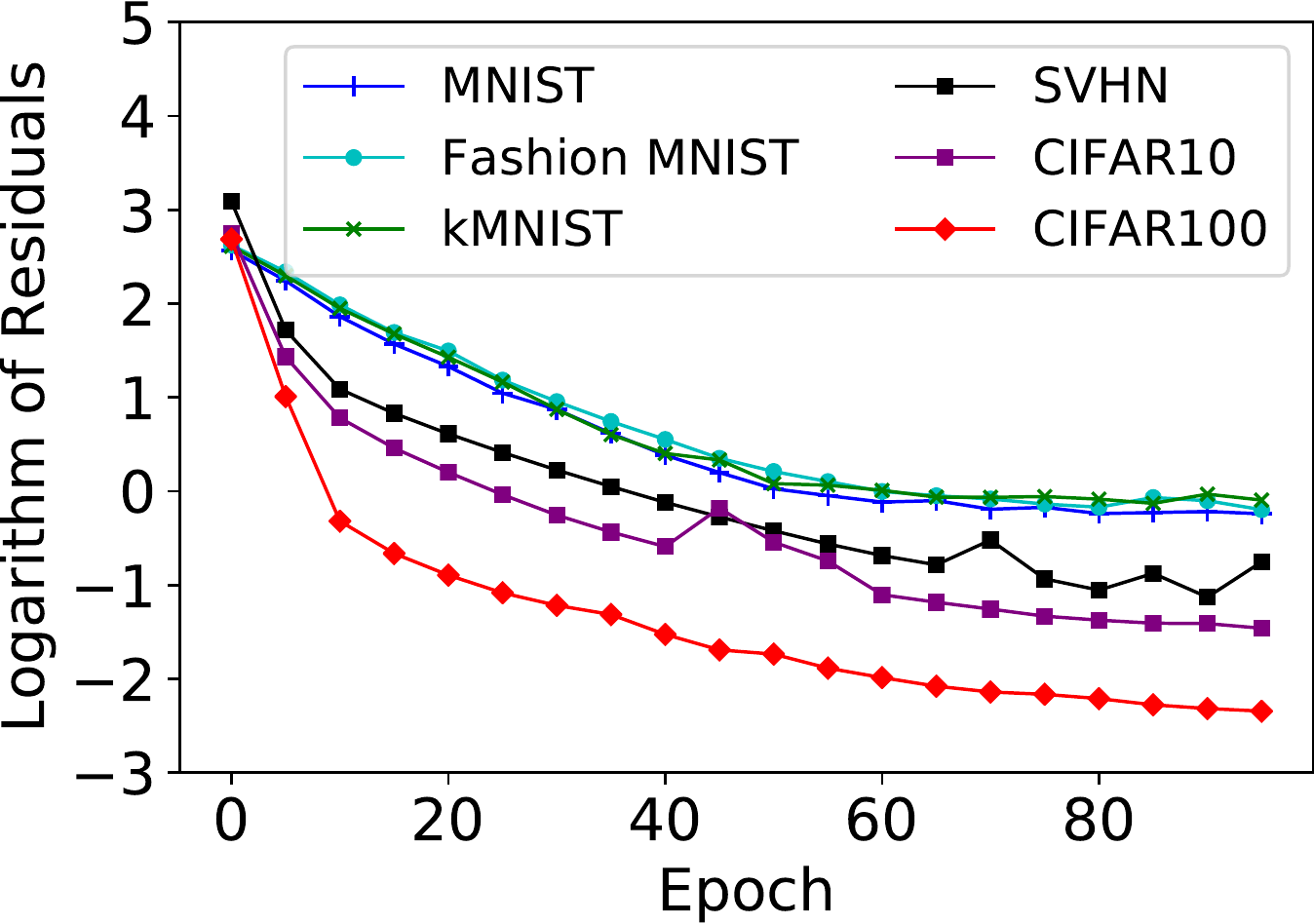}
     \centerline{(b). Residual versus epoch}
      \end{minipage}
      \caption{ The convergence of the proposed pdADMM: the objective decreases monotonously, and the residual converges to 0.}
      \label{fig:convergence}
      \vspace{-0.3cm}
 \end{figure*}
 To validate the convergence of the proposed pdADMM, we set up a feed-forward neural network with 9 hidden layers, each of which has 500 neurons. The Rectified linear unit (ReLU) was used for the activation function for both network structures. The loss function was set as the cross-entropy loss. The number of epoch was set to $100$. $\nu$ and $\rho$ were both set to $0.1$. As shown in Figure \ref{fig:convergence}, the objective keeps decreasing monotonically on all six datasets, and the residual converges sublinearly to 0, which are consistent with Theorems \ref{theo: global convergence} and \ref{theo: convergence rate}.
 \subsection{Accuracy}
 \begin{figure*}
    \begin{minipage}{0.3\linewidth}
    \includegraphics[width=\linewidth]{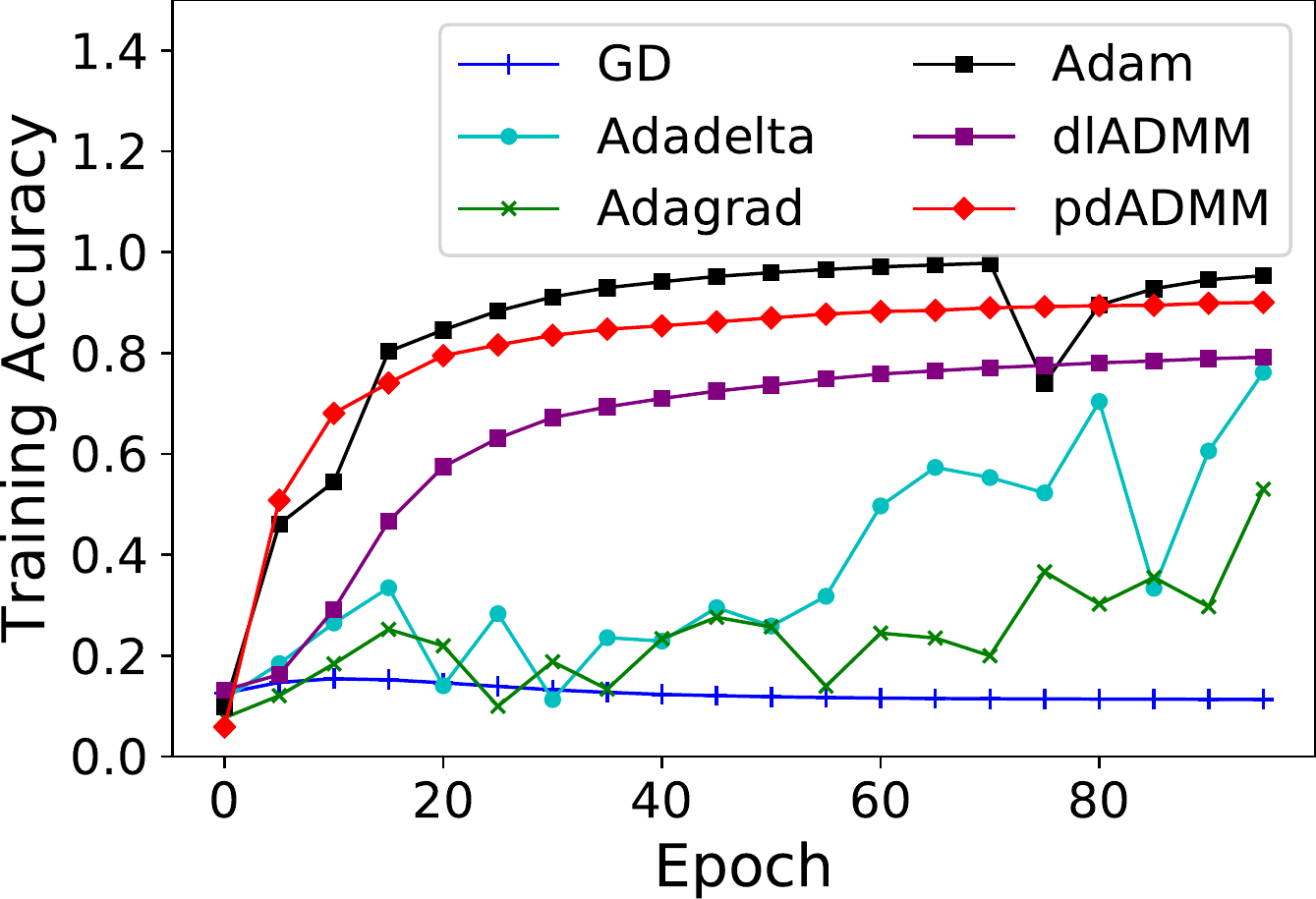}
    \centerline{(a). MNIST}
    \end{minipage}
    \hfill
    \begin{minipage}{0.3\linewidth}
    \includegraphics[width=\columnwidth]{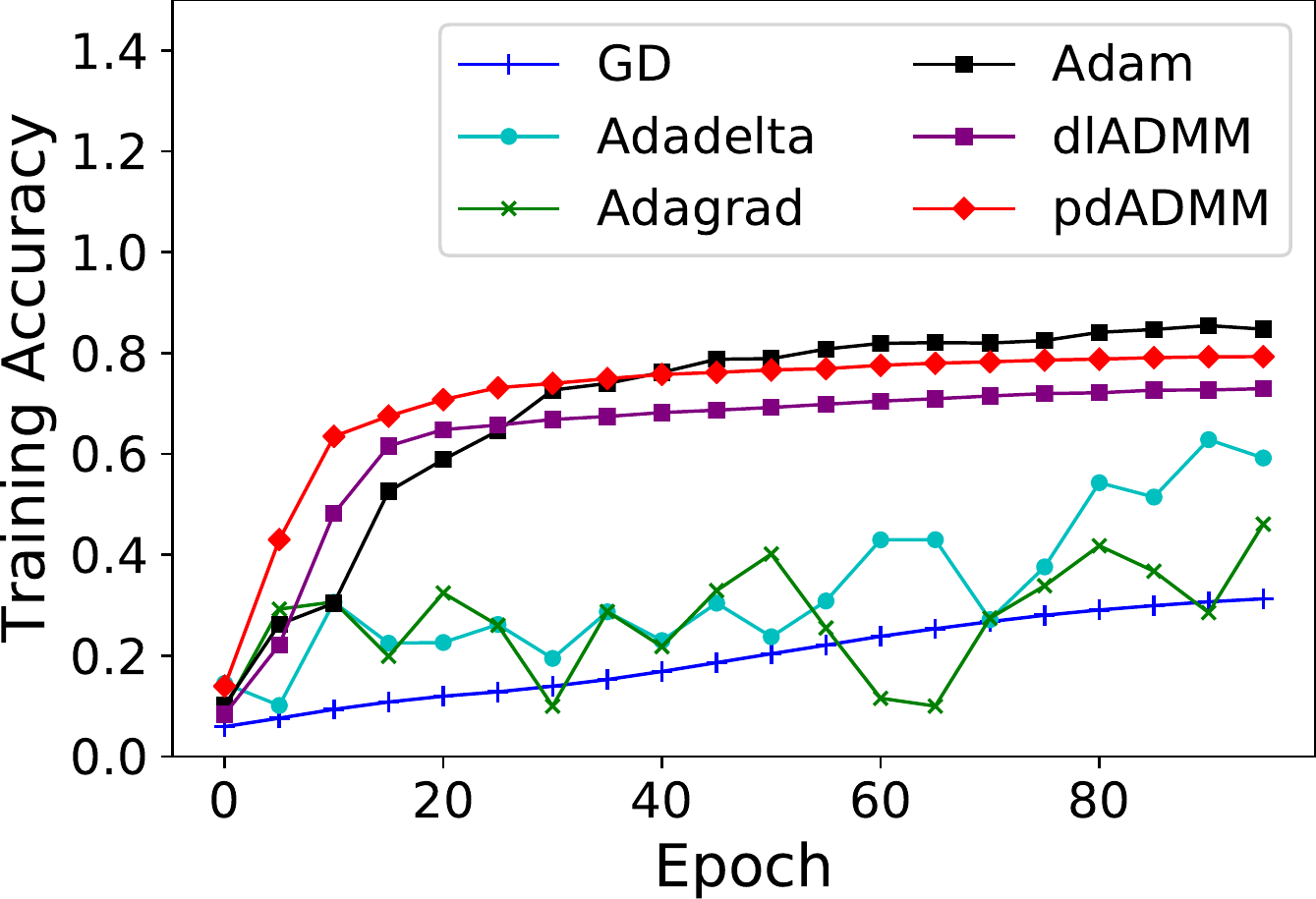}
    \centerline{(b). Fashion MNIST}
    \end{minipage}
    \hfill
    \begin{minipage}{0.3\linewidth}
    \includegraphics[width=\linewidth]{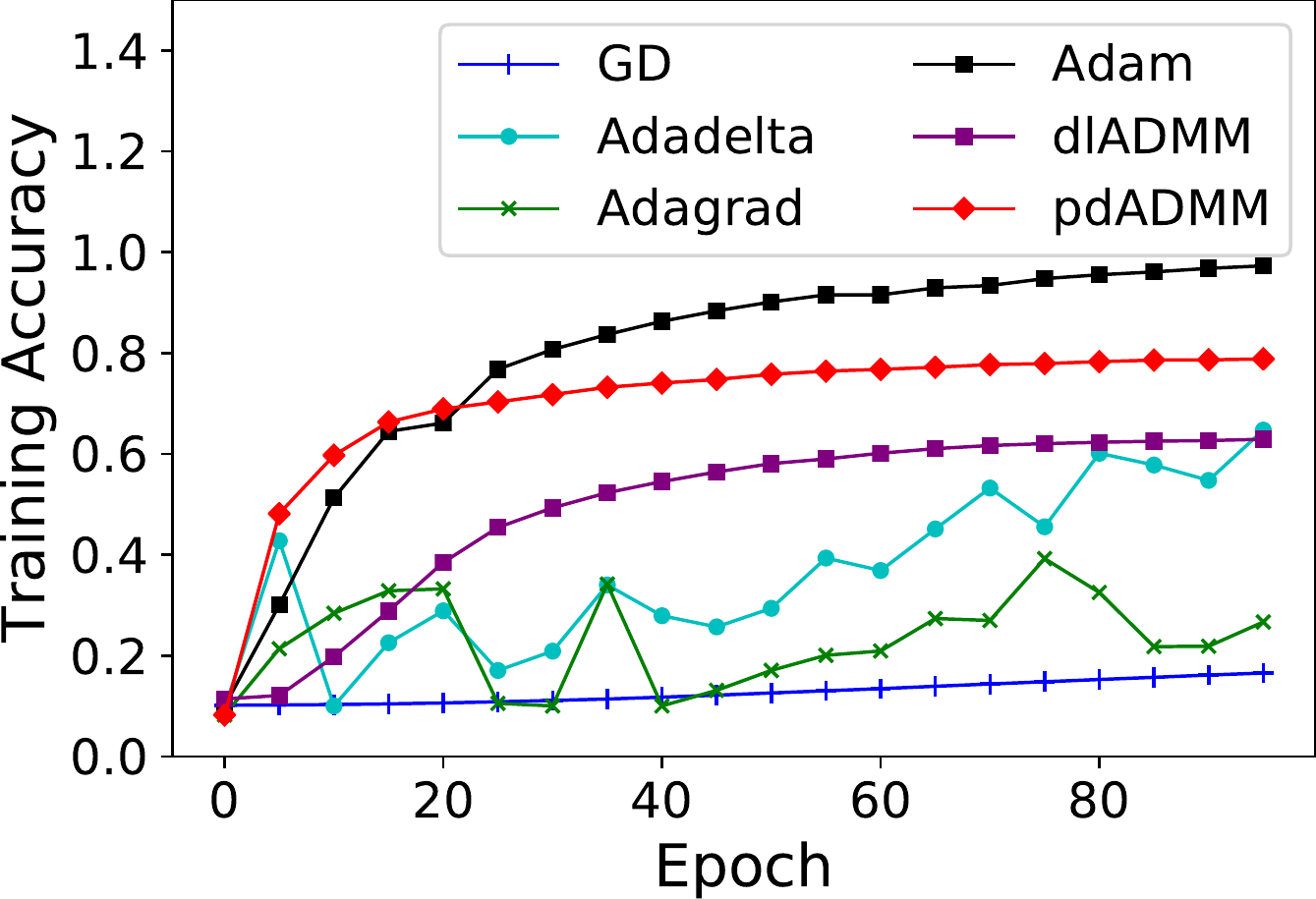}
        \centerline{(c). kMNIST}
    \end{minipage}
    \vfill
    \begin{minipage}{0.3\linewidth}
    \includegraphics[width=\linewidth]{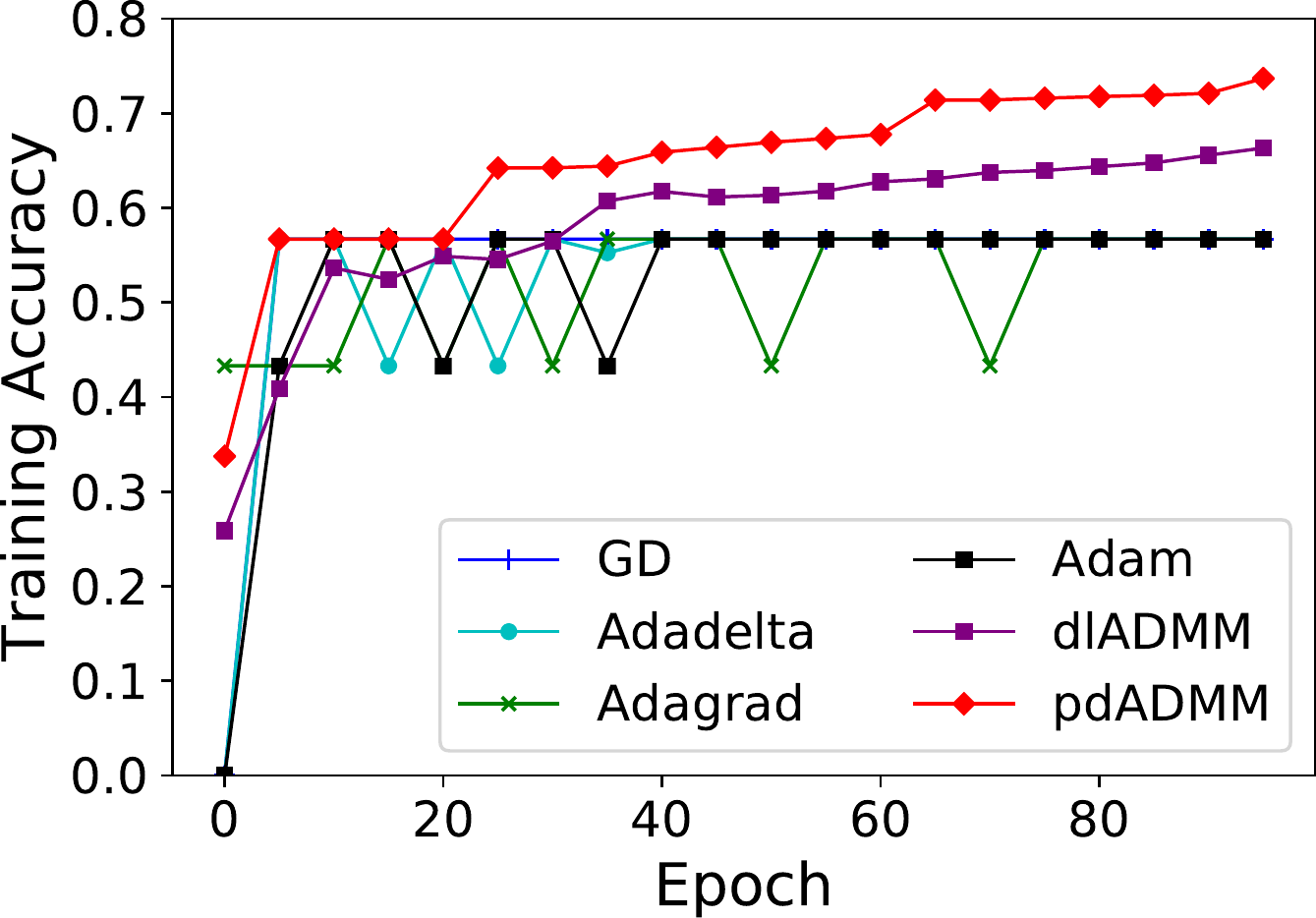}
    \centerline{(d). SVHN}
    \end{minipage}
    \hfill
    \begin{minipage}{0.3\linewidth}
    \includegraphics[width=\linewidth]{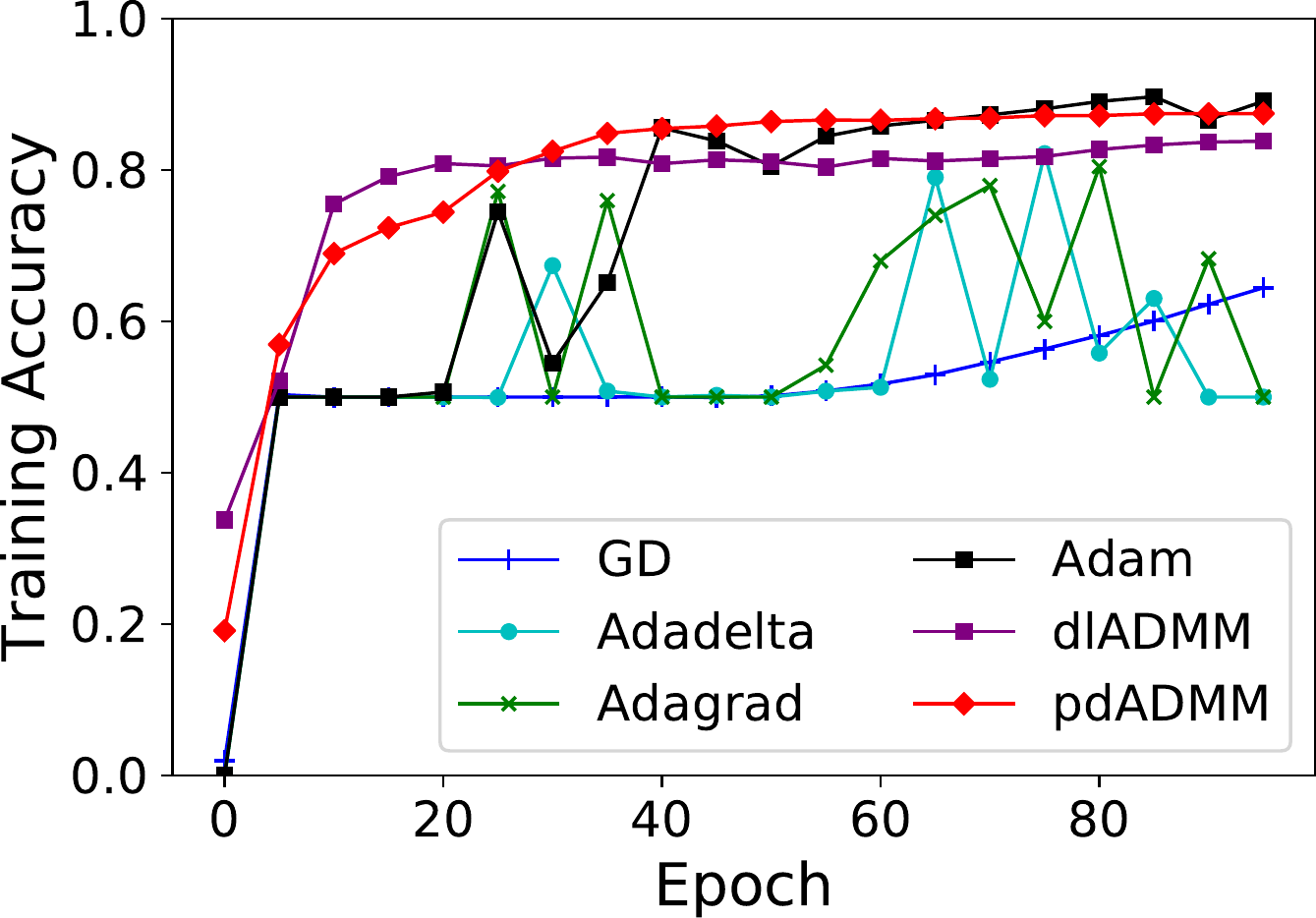}
        \centerline{(e). CIFAR10}
    \end{minipage}
    \hfill
    \begin{minipage}{0.3\linewidth}
    \includegraphics[width=\linewidth]{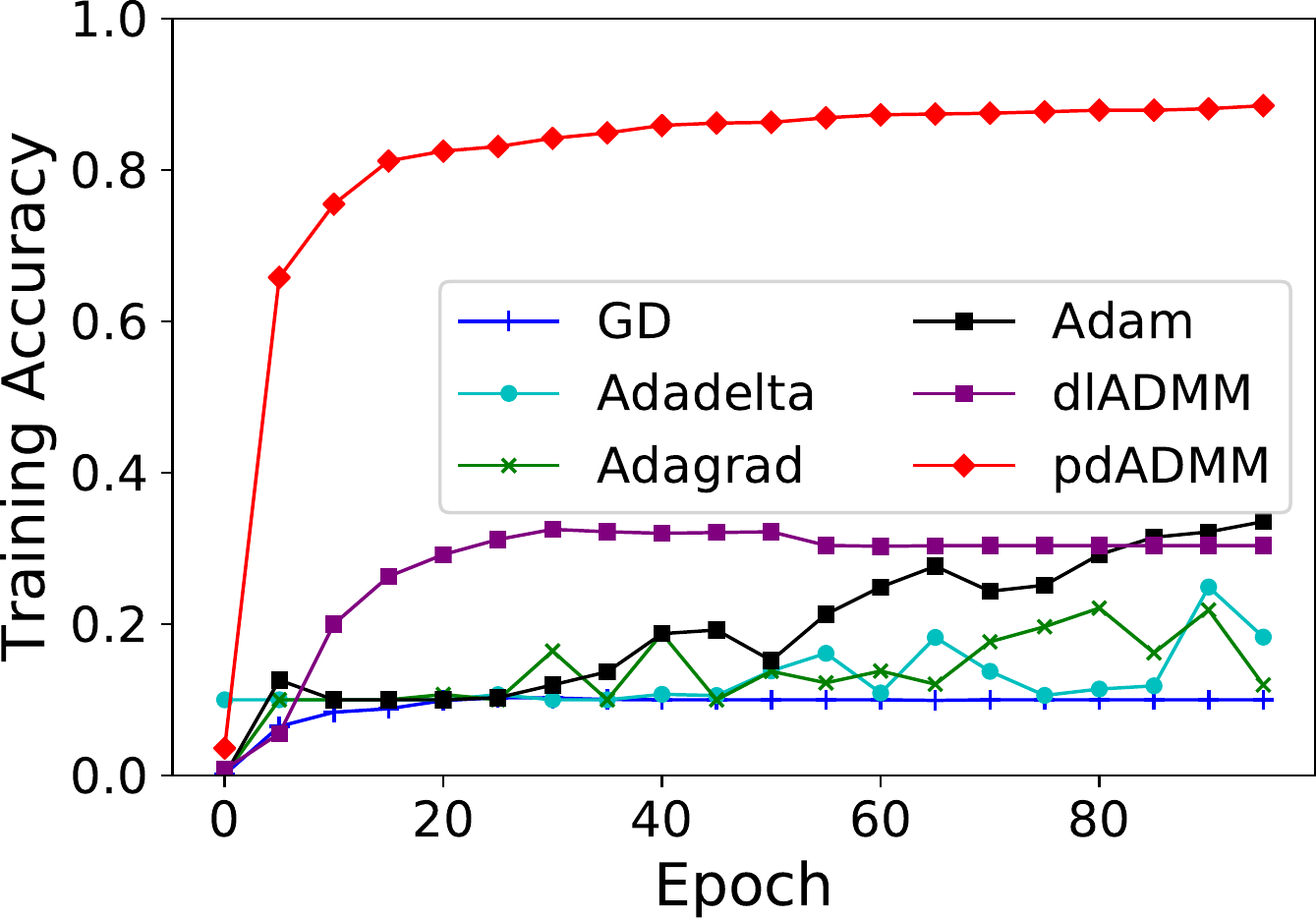}
            \centerline{(f). CIFAR100}
    \end{minipage}
    \caption{Training Accuracy of all methods: pdADMM outperformed most comparison methods; SGD-type methods suffer from gradient vanishing problems.}
    \label{fig:train performance}
    \vspace{-0.5cm}
\end{figure*}
\begin{figure*}
    \begin{minipage}{0.3\linewidth}
    \includegraphics[width=\linewidth]{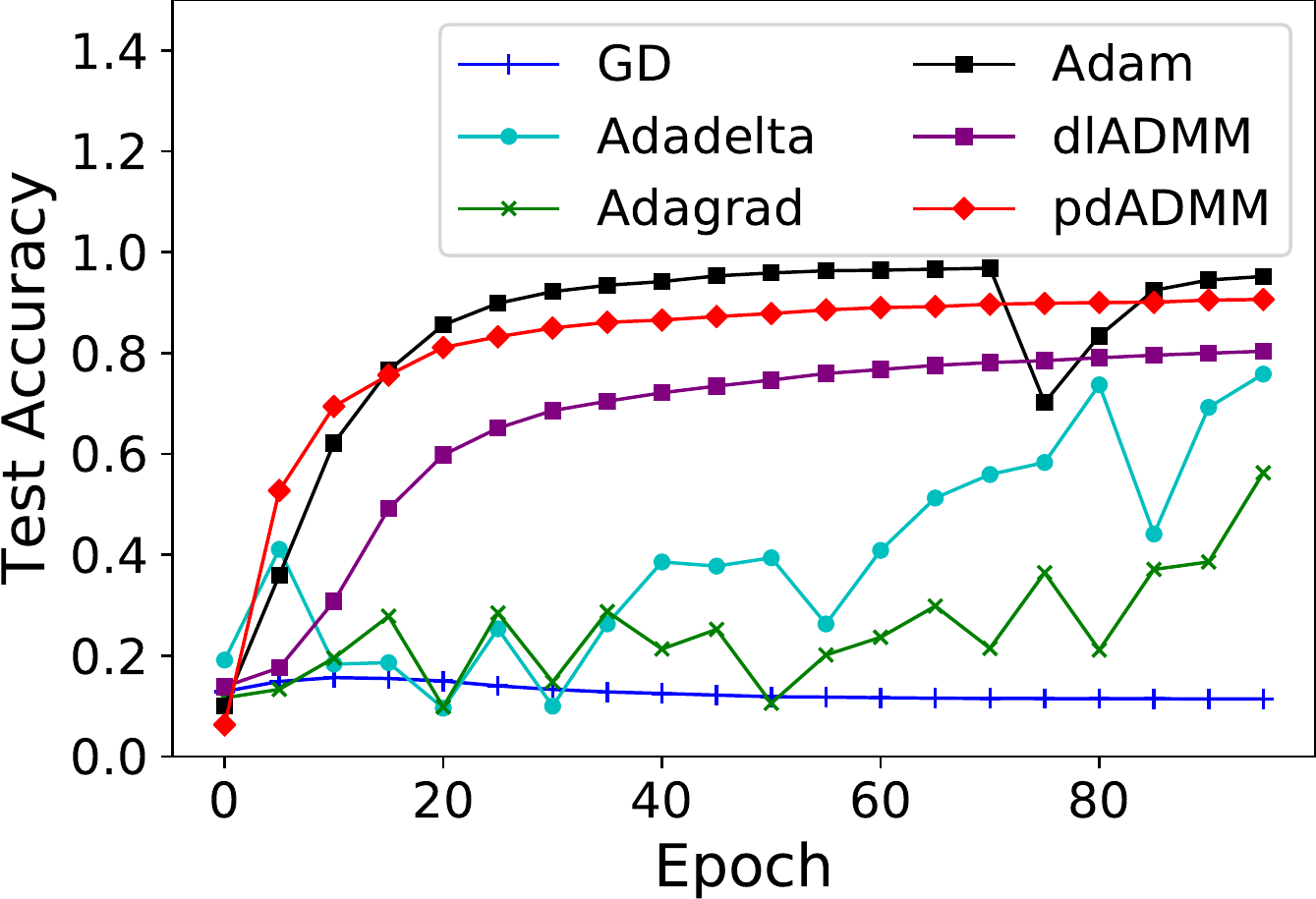}
    \centerline{(a). MNIST}
    \end{minipage}
    \hfill
    \begin{minipage}{0.3\linewidth}
    \includegraphics[width=\columnwidth]{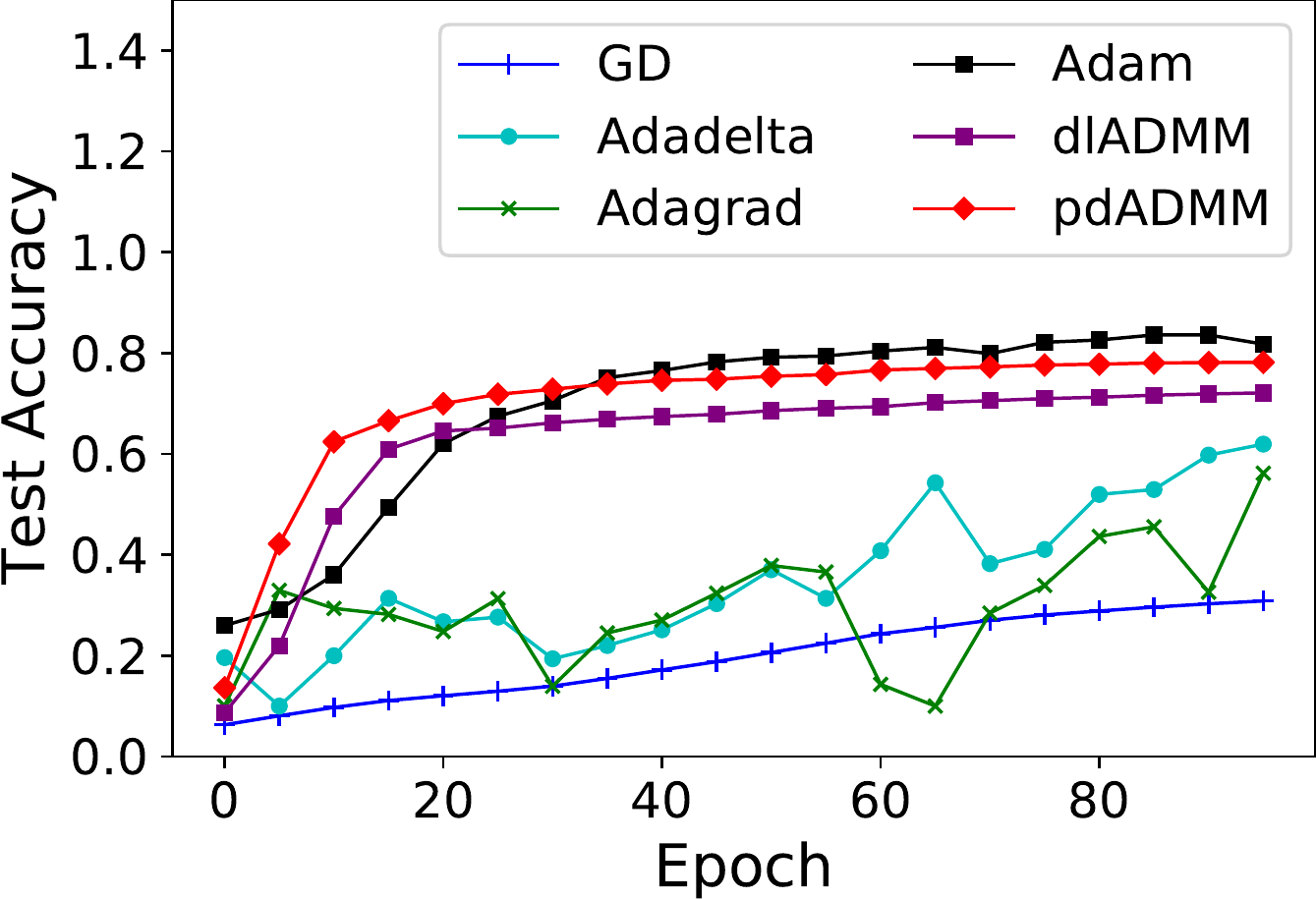}
    \centerline{(b). Fashion MNIST}
    \end{minipage}
    \hfill
    \begin{minipage}{0.3\linewidth}
    \includegraphics[width=\linewidth]{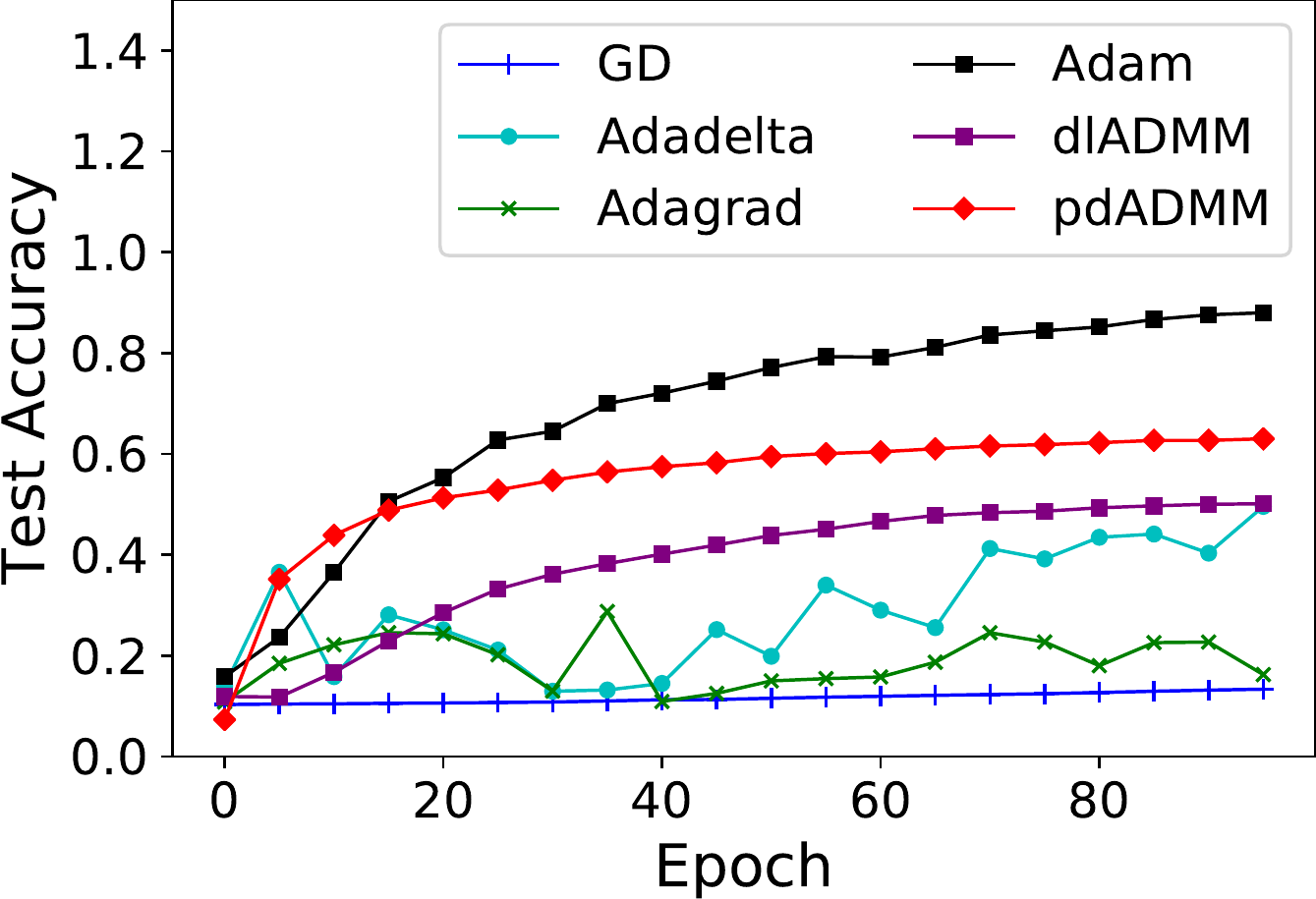}
        \centerline{(c). kMNIST}
    \end{minipage}
    \vfill
    \begin{minipage}{0.3\linewidth}
    \includegraphics[width=\linewidth]{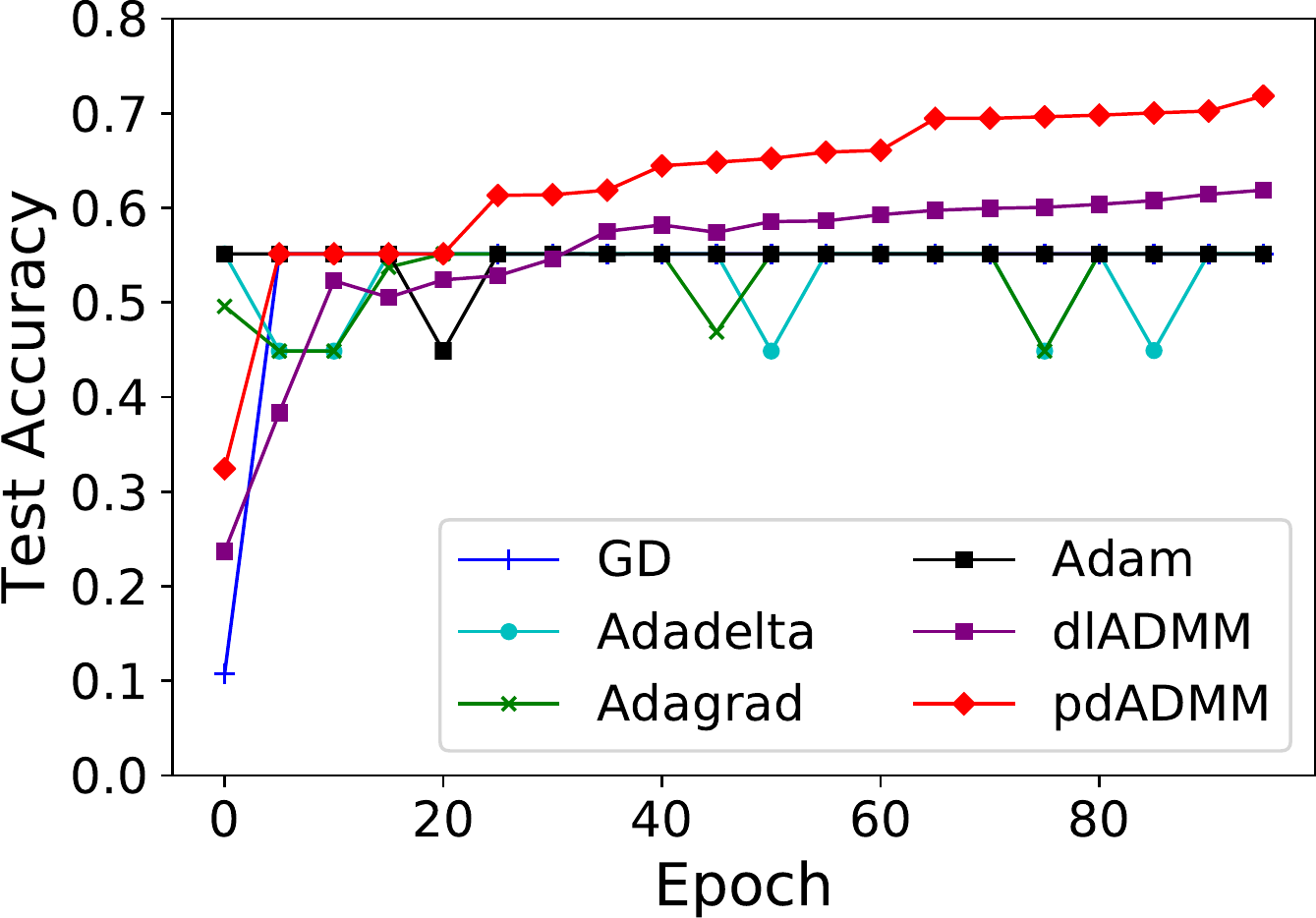}
    \centerline{(d). SVHN}
    \end{minipage}
    \hfill
    \begin{minipage}{0.3\linewidth}
    \includegraphics[width=\linewidth]{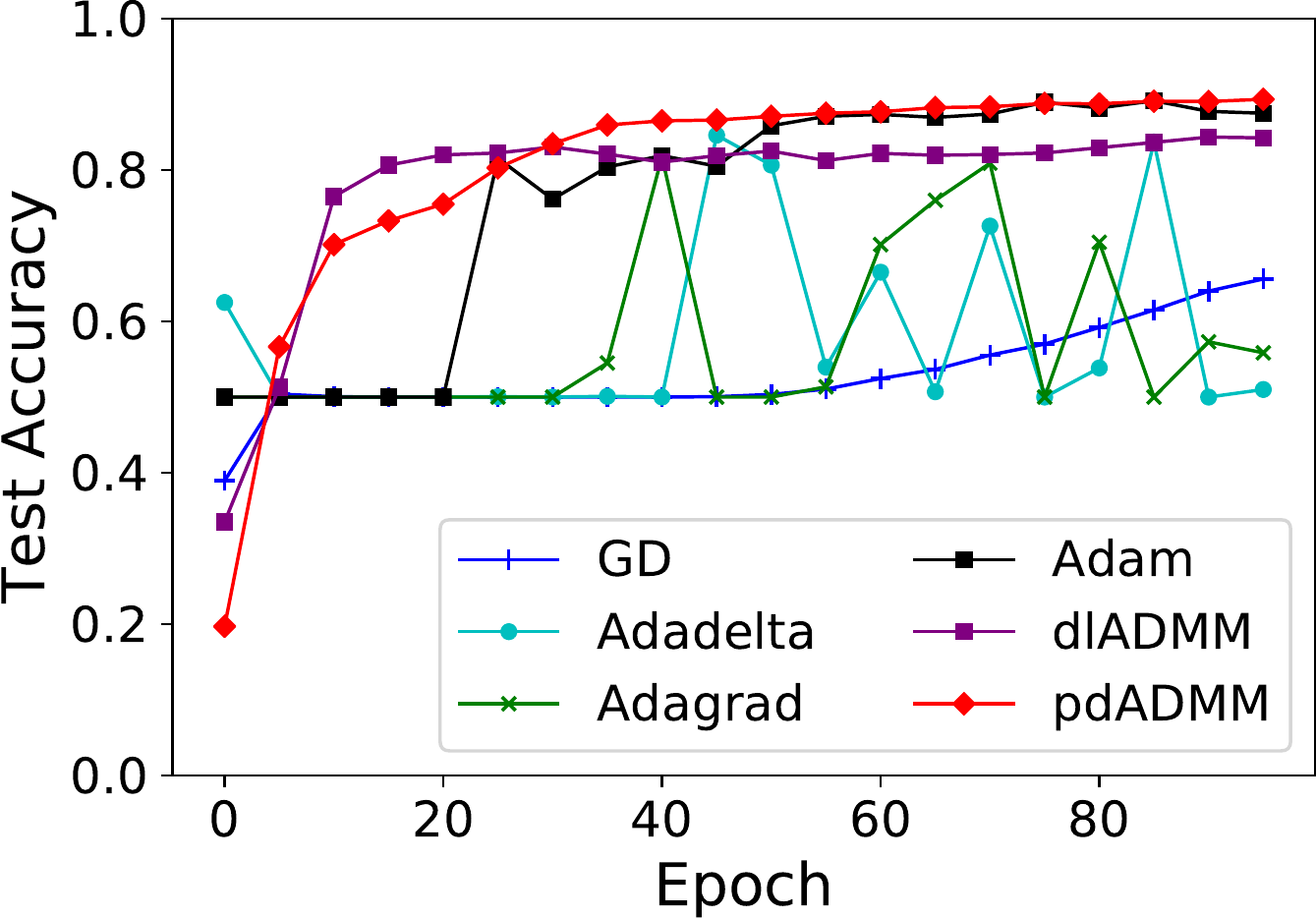}
        \centerline{(e). CIFAR10}
    \end{minipage}
    \hfill
    \begin{minipage}{0.3\linewidth}
    \includegraphics[width=\linewidth]{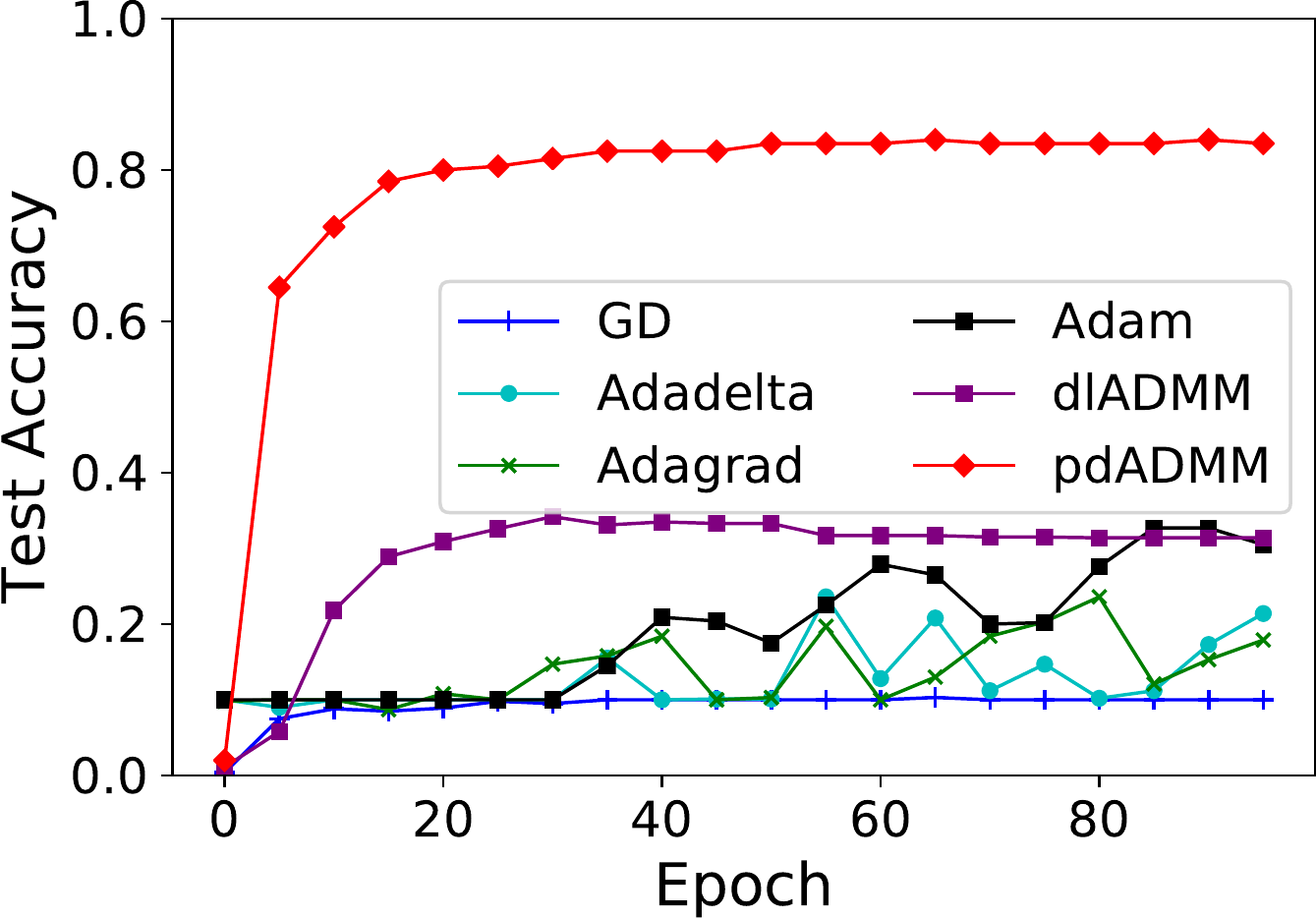}
            \centerline{(f). CIFAR100}
    \end{minipage}
    \caption{Test Accuracy of all methods: pdADMM outperformed most comparison methods.}
    \label{fig:test performance}
    \vspace{-0.7cm}
\end{figure*}

\subsubsection{Experimental Settings}
\indent In order to evaluate accuracy, we used the same architecture as the previous section.  $\nu$ and $\rho$ were set to $10^{-4}$ in order to maximize the performance of training data. The number of epoch was set to $100$. In this experiment, the full batch was used for training. As suggested by the training strategy in Section \ref{sec:subproblem}, we firstly trained a feed-forward neural network with five hidden layers, and then all layers were involved in training.
\subsubsection{Comparison Methods}
\indent SGD and its variants are state-of-the-art methods and hence were served as comparison methods. For SGD-based methods, the full batch dataset is used for training models. All parameters were chosen by maximizing the accuracy of training datasets. The baselines are described as follows: \\
\indent 1. Gradient Descent (GD) \cite{bottou2010large}. The GD and its variants are the most popular deep learning optimizers, whose convergence has been studied extensively in the literature. The learning rate of GD was set to $0.01$.\\
\indent 2. Adaptive learning rate method (Adadelta) \cite{zeiler2012adadelta}. The Adadelta is proposed to overcome the sensitivity to hyperparameter selection. The learning rate of Adadelta was set to $1$.\\
\indent 3. Adaptive gradient algorithm (Adagrad) \cite{duchi2011adaptive}. Adagrad is an improved version of SGD: rather than fixing the learning rate during training, it adapts the learning rate to the hyperparameter. The learning rate of Adagrad was set to $0.1$.\\
\indent 4. Adaptive momentum estimation (Adam) \cite{kingma2014adam}. Adam is the most popular optimization method for deep learning models. It estimates the first and second momentum to correct the biased gradient and thus makes convergence fast. The learning rate of Adam was set to $0.001$.\\
\indent 5. Deep learning Alternaing Direction Method of Multipliers (dlADMM) \cite{wang2019admm}. The dlADMM is an improvement of the previous ADMM implementation \cite{taylor2016training}. It is guaranteed to converge to a critical point with a rate of $o(1/k)$. $\rho$ and $\nu$ were both set to $10^{-6}$.
\subsubsection{Performance}
\indent In this section, the performance of the proposed pdADMM is analyzed against comparison methods. Figures \ref{fig:train performance} and \ref{fig:test performance} show the training and test accuracy of the proposed pdADMM against comparison methods on six datasets, respectively. X-axis and Y-axis represent epoch and training accuracy, respectively. Overall, the pdADMM outperformed most of comparison methods: it performed the best on the SVHN, CIFAR10 and CIFAR100 datasets, while was only secondary to Adam on the MNIST, Fashion MNIST and kMNIST datasets. The performance gap is particularly obvious on the CIFAR100 dataset. Most SGD-type methods suffered from gradient vanishing in the deep neural network, and struggled to find an optimum: for example, GD can only reach $10\%$ training accuracy on the  MNIST dataset; Adagrad and Adadelta performed better then GD, but took many epochs to escape saddle points. The dlADMM can avoid gradient vanishing problem, however, it performed worse than the proposed pdADMM on all datasets. Adam performed the best on three MNIST-like datasets, but performed worse than the proposed pdADMM on three other datasets, which are hard to train with high accuracy than three MNIST-like datasets. This indicates that the proposed pdADMM may be more suitable for training hard datasets than Adam.
\vspace{-0.2cm}
\section{Conclusion}
\label{sec:conclusion}
 Alternating Direction Method of Multipliers (ADMM) is considered to be a good alternative to Stochastic gradient descent (SGD) for training deep neural networks. In this paper, we propose a novel parallel deep learning Alternating Direction Method of Multipliers (pdADMM) to achieve layer parallelism.  The proposed pdADMM is guaranteed to converge to a critical solution under mild conditions. Experiments on benchmark datasets demonstrate that our proposed pdADMM can lead to a huge speedup when training a deep feed-forward neural network, it also outperformed others on six benchmark datasets.
 \vspace{-0.3cm}
 \section*{Acknowledgement}
 \vspace{-0.2cm}
\indent This work was supported by the National Science Foundation (NSF) Grant No. 1755850, No. 1841520, No. 2007716, No. 2007976, No. 1942594, No. 1907805, a Jeffress Memorial Trust Award, NVIDIA GPU Grant, and Design Knowledge Company (subcontract number: 10827.002.120.04).
 \small
\bibliography{example_paper}

\begin{thebibliography}{10}
\providecommand{\url}[1]{#1}
\csname url@samestyle\endcsname
\providecommand{\newblock}{\relax}
\providecommand{\bibinfo}[2]{#2}
\providecommand{\BIBentrySTDinterwordspacing}{\spaceskip=0pt\relax}
\providecommand{\BIBentryALTinterwordstretchfactor}{4}
\providecommand{\BIBentryALTinterwordspacing}{\spaceskip=\fontdimen2\font plus
\BIBentryALTinterwordstretchfactor\fontdimen3\font minus
  \fontdimen4\font\relax}
\providecommand{\BIBforeignlanguage}[2]{{%
\expandafter\ifx\csname l@#1\endcsname\relax
\typeout{** WARNING: IEEEtran.bst: No hyphenation pattern has been}%
\typeout{** loaded for the language `#1'. Using the pattern for}%
\typeout{** the default language instead.}%
\else
\language=\csname l@#1\endcsname
\fi
#2}}
\providecommand{\BIBdecl}{\relax}
\BIBdecl

\bibitem{taylor2016training}
G.~Taylor, R.~Burmeister, Z.~Xu, B.~Singh, A.~Patel, and T.~Goldstein,
  ``Training neural networks without gradients: A scalable admm approach,'' in
  \emph{International conference on machine learning}, 2016, pp. 2722--2731.

\bibitem{wang2019admm}
J.~Wang, F.~Yu, X.~Chen, and L.~Zhao, ``Admm for efficient deep learning with
  global convergence,'' in \emph{Proceedings of the 25th ACM SIGKDD
  International Conference on Knowledge Discovery \& Data Mining}, 2019, pp.
  111--119.

\bibitem{boyd2011distributed}
S.~Boyd, N.~Parikh, E.~Chu, B.~Peleato, J.~Eckstein \emph{et~al.},
  ``Distributed optimization and statistical learning via the alternating
  direction method of multipliers,'' \emph{Foundations and
  Trends{\textregistered} in Machine learning}, vol.~3, no.~1, pp. 1--122,
  2011.

\bibitem{mota2012distributed}
J.~F. Mota, J.~M. Xavier, P.~M. Aguiar, and M.~P{\"u}schel, ``Distributed admm
  for model predictive control and congestion control,'' in \emph{2012 IEEE
  51st IEEE Conference on Decision and Control (CDC)}.\hskip 1em plus 0.5em
  minus 0.4em\relax IEEE, 2012, pp. 5110--5115.

\bibitem{chang2016proximal}
T.-H. Chang, ``A proximal dual consensus admm method for multi-agent
  constrained optimization,'' \emph{IEEE Transactions on Signal Processing},
  vol.~64, no.~14, pp. 3719--3734, 2016.

\bibitem{chang2014multi}
T.-H. Chang, M.~Hong, and X.~Wang, ``Multi-agent distributed optimization via
  inexact consensus admm,'' \emph{IEEE Transactions on Signal Processing},
  vol.~63, no.~2, pp. 482--497, 2014.

\bibitem{shi2014linear}
W.~Shi, Q.~Ling, K.~Yuan, G.~Wu, and W.~Yin, ``On the linear convergence of the
  admm in decentralized consensus optimization,'' \emph{IEEE Transactions on
  Signal Processing}, vol.~62, no.~7, pp. 1750--1761, 2014.

\bibitem{xu2017adaptive}
Z.~Xu, G.~Taylor, H.~Li, M.~A. Figueiredo, X.~Yuan, and T.~Goldstein,
  ``Adaptive consensus admm for distributed optimization,'' in
  \emph{Proceedings of the 34th International Conference on Machine
  Learning-Volume 70}.\hskip 1em plus 0.5em minus 0.4em\relax JMLR. org, 2017,
  pp. 3841--3850.

\bibitem{zhu2016quantized}
S.~Zhu, M.~Hong, and B.~Chen, ``Quantized consensus admm for multi-agent
  distributed optimization,'' in \emph{2016 IEEE International Conference on
  Acoustics, Speech and Signal Processing (ICASSP)}.\hskip 1em plus 0.5em minus
  0.4em\relax IEEE, 2016, pp. 4134--4138.

\bibitem{zhang2014asynchronous}
R.~Zhang and J.~Kwok, ``Asynchronous distributed admm for consensus
  optimization,'' in \emph{International Conference on Machine Learning}, 2014,
  pp. 1701--1709.

\bibitem{wei2012distributed}
E.~Wei and A.~Ozdaglar, ``Distributed alternating direction method of
  multipliers,'' in \emph{2012 IEEE 51st IEEE Conference on Decision and
  Control (CDC)}.\hskip 1em plus 0.5em minus 0.4em\relax IEEE, 2012, pp.
  5445--5450.

\bibitem{chang2016asynchronous1}
T.-H. Chang, M.~Hong, W.-C. Liao, and X.~Wang, ``Asynchronous distributed admm
  for large-scale optimization—part i: Algorithm and convergence analysis,''
  \emph{IEEE Transactions on Signal Processing}, vol.~64, no.~12, pp.
  3118--3130, 2016.

\bibitem{chang2016asynchronous2}
T.-H. Chang, W.-C. Liao, M.~Hong, and X.~Wang, ``Asynchronous distributed admm
  for large-scale optimization—part ii: Linear convergence analysis and
  numerical performance,'' \emph{IEEE Transactions on Signal Processing},
  vol.~64, no.~12, pp. 3131--3144, 2016.

\bibitem{hong2014distributed}
M.~Hong, ``A distributed, asynchronous, and incremental algorithm for nonconvex
  optimization: an admm approach,'' \emph{IEEE Transactions on Control of
  Network Systems}, vol.~5, no.~3, pp. 935--945, 2017.

\bibitem{kumar2016asynchronous}
S.~Kumar, R.~Jain, and K.~Rajawat, ``Asynchronous optimization over
  heterogeneous networks via consensus admm,'' \emph{IEEE Transactions on
  Signal and Information Processing over Networks}, vol.~3, no.~1, pp.
  114--129, 2016.

\bibitem{magnusson2015convergence}
S.~Magn{\'u}sson, P.~C. Weeraddana, M.~G. Rabbat, and C.~Fischione, ``On the
  convergence of alternating direction lagrangian methods for nonconvex
  structured optimization problems,'' \emph{IEEE Transactions on Control of
  Network Systems}, vol.~3, no.~3, pp. 296--309, 2015.

\bibitem{li2015global}
G.~Li and T.~K. Pong, ``Global convergence of splitting methods for nonconvex
  composite optimization,'' \emph{SIAM Journal on Optimization}, vol.~25,
  no.~4, pp. 2434--2460, 2015.

\bibitem{wang2015global}
Y.~Wang, W.~Yin, and J.~Zeng, ``Global convergence of admm in nonconvex
  nonsmooth optimization,'' \emph{Journal of Scientific Computing}, pp. 1--35,
  2015.

\bibitem{hong2016convergence}
M.~Hong, Z.-Q. Luo, and M.~Razaviyayn, ``Convergence analysis of alternating
  direction method of multipliers for a family of nonconvex problems,''
  \emph{SIAM Journal on Optimization}, vol.~26, no.~1, pp. 337--364, 2016.

\bibitem{wang2019multi}
J.~Wang, L.~Zhao, and L.~Wu, ``Multi-convex inequality-constrained alternating
  direction method of multipliers,'' \emph{arXiv preprint arXiv:1902.10882},
  2019.

\bibitem{liu2019linearized}
Q.~Liu, X.~Shen, and Y.~Gu, ``Linearized admm for nonconvex nonsmooth
  optimization with convergence analysis,'' \emph{IEEE Access}, vol.~7, pp.
  76\,131--76\,144, 2019.

\bibitem{wang2017nonconvex}
J.~Wang and L.~Zhao, ``Nonconvex generalization of admm for nonlinear equality
  constrained problems,'' \emph{arXiv preprint arXiv:1705.03412}, 2017.

\bibitem{xie2019differentiable}
X.~Xie, J.~Wu, G.~Liu, Z.~Zhong, and Z.~Lin, ``Differentiable linearized
  admm,'' in \emph{International Conference on Machine Learning}, 2019, pp.
  6902--6911.

\bibitem{wang2019accelerated}
J.~Wang, F.~Yu, X.~Chen, and L.~Zhao, ``Accelerated gradient-free neural
  network training by multi-convex alternating optimization,'' 2019.

\bibitem{chartrand2013nonconvex}
R.~Chartrand and B.~Wohlberg, ``A nonconvex admm algorithm for group sparsity
  with sparse groups,'' in \emph{Acoustics, Speech and Signal Processing
  (ICASSP), 2013 IEEE International Conference on}.\hskip 1em plus 0.5em minus
  0.4em\relax IEEE, 2013, pp. 6009--6013.

\bibitem{hajinezhad2015nonconvex}
D.~Hajinezhad and M.~Hong, ``Nonconvex alternating direction method of
  multipliers for distributed sparse principal component analysis,'' in
  \emph{2015 IEEE Global Conference on Signal and Information Processing
  (GlobalSIP)}.\hskip 1em plus 0.5em minus 0.4em\relax IEEE, 2015, pp.
  255--259.

\bibitem{guo2017convergence}
K.~Guo, D.~Han, D.~Z. Wang, and T.~Wu, ``Convergence of admm for multi-block
  nonconvex separable optimization models,'' \emph{Frontiers of Mathematics in
  China}, vol.~12, no.~5, pp. 1139--1162, 2017.

\bibitem{themelis2020douglas}
A.~Themelis and P.~Patrinos, ``Douglas--rachford splitting and admm for
  nonconvex optimization: Tight convergence results,'' \emph{SIAM Journal on
  Optimization}, vol.~30, no.~1, pp. 149--181, 2020.

\bibitem{wang2020tssm}
J.~Wang, Z.~Chai, Y.~Chen, and L.~Zhao, ``Tunable subnetwork splitting for
  model-parallelism of neural network training,'' in \emph{ICML 2020 Workshop:
  Beyond First Order Methods in Machine Learning}, 2020.

\bibitem{wang2019opt}
J.~Wang and L.~Zhao, ``The application of multi-block admm on isotonic
  regression problems,'' \emph{arXiv preprint arXiv:1903.01054}, 2019.

\bibitem{wen2017terngrad}
W.~Wen, C.~Xu, F.~Yan, C.~Wu, Y.~Wang, Y.~Chen, and H.~Li, ``Terngrad: Ternary
  gradients to reduce communication in distributed deep learning,'' in
  \emph{Advances in neural information processing systems}, 2017, pp.
  1509--1519.

\bibitem{sergeev2018horovod}
A.~Sergeev and M.~Del~Balso, ``Horovod: fast and easy distributed deep learning
  in tensorflow,'' 2018.

\bibitem{ooi2015singa}
B.~C. Ooi, K.-L. Tan, S.~Wang, W.~Wang, Q.~Cai, G.~Chen, J.~Gao, Z.~Luo, A.~K.
  Tung, Y.~Wang \emph{et~al.}, ``Singa: A distributed deep learning platform,''
  in \emph{Proceedings of the 23rd ACM international conference on
  Multimedia}.\hskip 1em plus 0.5em minus 0.4em\relax ACM, 2015, pp. 685--688.

\bibitem{chen2015mxnet}
T.~Chen, M.~Li, Y.~Li, M.~Lin, N.~Wang, M.~Wang, T.~Xiao, B.~Xu, C.~Zhang, and
  Z.~Zhang, ``Mxnet: A flexible and efficient machine learning library for
  heterogeneous distributed systems,'' 2015.

\bibitem{hashemitictac}
S.~H. Hashemi, S.~A. Jyothi, and R.~H. Campbell, ``Tictac: Accelerating
  distributed deep learning with communication scheduling,'' in
  \emph{Proceedings of the 2nd SysML Conference}, 2019.

\bibitem{zhang2017poseidon}
H.~Zhang, Z.~Zheng, S.~Xu, W.~Dai, Q.~Ho, X.~Liang, Z.~Hu, J.~Wei, P.~Xie, and
  E.~P. Xing, ``Poseidon: An efficient communication architecture for
  distributed deep learning on $\{$GPU$\}$ clusters,'' in \emph{2017
  $\{$USENIX$\}$ Annual Technical Conference ($\{$USENIX$\}$$\{$ATC$\}$ 17)},
  2017, pp. 181--193.

\bibitem{zinkevich2010parallelized}
M.~Zinkevich, M.~Weimer, L.~Li, and A.~J. Smola, ``Parallelized stochastic
  gradient descent,'' in \emph{Advances in neural information processing
  systems}, 2010, pp. 2595--2603.

\bibitem{parpas2019predict}
P.~Parpas and C.~Muir, ``Predict globally, correct locally: Parallel-in-time
  optimal control of neural networks,'' 2019.

\bibitem{huo2018training}
Z.~Huo, B.~Gu, and H.~Huang, ``Training neural networks using features
  replay,'' in \emph{Advances in Neural Information Processing Systems}, 2018,
  pp. 6659--6668.

\bibitem{zhuang2019fully}
H.~Zhuang, Y.~Wang, Q.~Liu, and Z.~Lin, ``Fully decoupled neural network
  learning using delayed gradients,'' 2019.

\bibitem{beck2009fast}
A.~Beck and M.~Teboulle, ``A fast iterative shrinkage-thresholding algorithm
  for linear inverse problems,'' \emph{SIAM journal on imaging sciences},
  vol.~2, no.~1, pp. 183--202, 2009.

\bibitem{rockafellar2009variational}
R.~T. Rockafellar and R.~J.-B. Wets, \emph{Variational analysis}.\hskip 1em
  plus 0.5em minus 0.4em\relax Springer Science \& Business Media, 2009, vol.
  317.

\bibitem{lecun1998gradient}
Y.~LeCun, L.~Bottou, Y.~Bengio, P.~Haffner \emph{et~al.}, ``Gradient-based
  learning applied to document recognition,'' \emph{Proceedings of the IEEE},
  vol.~86, no.~11, pp. 2278--2324, 1998.

\bibitem{chollet2015keras}
F.~Chollet \emph{et~al.}, ``Keras,'' https://keras.io, 2015.

\bibitem{xiao2017fashion}
H.~Xiao, K.~Rasul, and R.~Vollgraf, ``Fashion-mnist: a novel image dataset for
  benchmarking machine learning algorithms,'' 2017.

\bibitem{clanuwat2018deep}
T.~Clanuwat, M.~Bober-Irizar, A.~Kitamoto, A.~Lamb, K.~Yamamoto, and D.~Ha,
  ``Deep learning for classical japanese literature,'' in \emph{NeurIPS 2018
  Workshop on Machine Learning for Creativity and Design}, 2018.

\bibitem{SVHN}
Y.~Netzer, T.~Wang, A.~Coates, A.~Bissacco, B.~Wu, and A.~Y. Ng, ``Reading
  digits in natural images with unsupervised feature learning,'' in \emph{NIPS
  Workshop on Deep Learning and Unsupervised Feature Learning 2011}, 2011.

\bibitem{krizhevsky2009learning}
A.~Krizhevsky \emph{et~al.}, ``Learning multiple layers of features from tiny
  images,'' Citeseer, Tech. Rep., 2009.

\bibitem{bottou2010large}
L.~Bottou, ``Large-scale machine learning with stochastic gradient descent,''
  in \emph{Proceedings of COMPSTAT'2010}.\hskip 1em plus 0.5em minus
  0.4em\relax Springer, 2010, pp. 177--186.

\bibitem{zeiler2012adadelta}
M.~D. Zeiler, ``Adadelta: an adaptive learning rate method,'' 2012.

\bibitem{duchi2011adaptive}
J.~Duchi, E.~Hazan, and Y.~Singer, ``Adaptive subgradient methods for online
  learning and stochastic optimization,'' \emph{Journal of Machine Learning
  Research}, vol.~12, no. Jul, pp. 2121--2159, 2011.

\bibitem{kingma2014adam}
D.~P. Kingma and J.~Ba, ``Adam: {A} method for stochastic optimization,'' in
  \emph{3rd International Conference on Learning Representations, {ICLR} 2015,
  San Diego, CA, USA, May 7-9, 2015, Conference Track Proceedings}, Y.~Bengio
  and Y.~LeCun, Eds., 2015.

\bibitem{deng2017parallel}
W.~Deng, M.-J. Lai, Z.~Peng, and W.~Yin, ``Parallel multi-block admm with o
  (1/k) convergence,'' \emph{Journal of Scientific Computing}, vol.~71, no.~2,
  pp. 712--736, 2017.

\end{thebibliography}
\bibliographystyle{IEEEtran}
\textbf{Proof of Theorem \ref{theo: convergence rate}}
\begin{proof}
To prove this theorem, we will first show that $c_k$ satisfies two conditions: (1). $c_k\geq c_{k+1}$. (2). $\sum\nolimits_{k=0}^\infty c_k$ is bounded. Specifically, first, we have
\begin{align*}
    c_k&\!=\!\min\nolimits_{0\leq i\leq k}\!(\sum\nolimits_{l=2}^L (\tau^{i+1}_l/2)\Vert p^{i+1}_l\!-\!p^i_l\Vert^2_2\\&+\sum\nolimits_{l\!=\!1}^{L}(\theta^{i\!+\!1}_l/2)\Vert  W^{i\!+\!1}_l\!-\!W^i_l\Vert^2_2+\sum\nolimits_{l\!=\!1}^{L}(\nu/2)\Vert  b^{i\!+\!1}_l\!-\!b^i_l\Vert^2_2\\&+\sum\nolimits_{l=1}^{L-1} C_1\Vert z^{i+1}_l-z^i_l\Vert^2_2+(\nu/2)\Vert z^{i+1}_L-z^i_L\Vert^2_2\\&+\sum\nolimits_{l=1}^{L-1}C_2\Vert q^{i+1}_l-q^i_l\Vert^2_2) \\&\geq\!\min\nolimits_{0\leq i\leq k+1}\!(\sum\nolimits_{l=2}^L (\tau^{i+1}_l/2)\Vert p^{i+1}_l\!-\!p^i_l\Vert^2_2\\&+\sum\nolimits_{l\!=\!1}^{L}(\theta^{i\!+\!1}_l/2)\Vert  W^{i\!+\!1}_l\!-\!W^i_l\Vert^2_2+\sum\nolimits_{l\!=\!1}^{L}(\nu/2)\Vert  b^{i\!+\!1}_l\!-\!b^i_l\Vert^2_2\\&+\sum\nolimits_{l=1}^{L-1} C_1\Vert z^{i+1}_l-z^i_l\Vert^2_2+(\nu/2)\Vert z^{i+1}_L-z^i_L\Vert^2_2\\&+\sum\nolimits_{l=1}^{L-1}C_2\Vert q^{i+1}_l-q^i_l\Vert^2_2)\\&= c_{k+1}
\end{align*}
Therefore $c_k$ satisfies the first condition. Second,
\begin{align*}
    &\sum\nolimits_{k\!=\!0}^\infty c_k=\sum\nolimits_{k=0}^\infty\min\nolimits_{0\leq i\leq k}\!(\sum\nolimits_{l=2}^L (\tau^{i+1}_l/2)\Vert p^{i+1}_l\!-\!p^i_l\Vert^2_2\\&+\sum\nolimits_{l\!=\!1}^{L}(\theta^{i\!+\!1}_l/2)\Vert  W^{i\!+\!1}_l\!-\!W^i_l\Vert^2_2+\sum\nolimits_{l\!=\!1}^{L}(\nu/2)\Vert  b^{i\!+\!1}_l\!-\!b^i_l\Vert^2_2\\&+\sum\nolimits_{l=1}^{L-1} C_1\Vert z^{i+1}_l-z^i_l\Vert^2_2+(\nu/2)\Vert z^{i+1}_L-z^i_L\Vert^2_2\\&+\sum\nolimits_{l=1}^{L-1}C_2\Vert q^{i+1}_l-q^i_l\Vert^2_2)\\&\leq \sum\nolimits_{k=0}^\infty(\sum\nolimits_{l=2}^L (\tau^{k+1}_l/2)\Vert p^{k+1}_l\!-\!p^k_l\Vert^2_2\\&+\sum\nolimits_{l\!=\!1}^{L}(\theta^{k\!+\!1}_l/2)\Vert  W^{k\!+\!1}_l\!-\!W^k_l\Vert^2_2+\sum\nolimits_{l\!=\!1}^{L}(\nu/2)\Vert  b^{k\!+\!1}_l\!-\!b^k_l\Vert^2_2\\&+\sum\nolimits_{l=1}^{L-1} C_1\Vert z^{k+1}_l-z^k_l\Vert^2_2+(\nu/2)\Vert z^{k+1}_L-z^k_L\Vert^2_2\\&+\sum\nolimits_{l=1}^{L-1}C_2\Vert q^{k+1}_l-q^k_l\Vert^2_2) \\&\leq L_\rho(\textbf{p}^0,\textbf{W}^0,\textbf{b}^0,\textbf{z}^0,\textbf{q}^0,\textbf{u}^0)-L_\rho(\textbf{p}^*,\textbf{W}^*,\textbf{b}^*,\textbf{z}^*,\textbf{q}^*,\textbf{u}^*)\\&\text{(Lemma \ref{lemma:objective decrease})}
\end{align*}
So $c_{k}$ satisfies the second condition. Finally, since we have proved the first two conditions and the third one $c_k \geq 0$ is obvious, the convergence rate of $o(1/k)$ is proven (Lemma 1.2 in \cite{deng2017parallel}). 
\end{proof}
\normalsize
\newpage
\onecolumn
\begin{appendix}
\textbf{Preliminary Results}
\begin{lemma}
\label{lemma: q_opt}
It holds for every $k\in \mathbb{N}$ and $l=1,\cdots,L-1$ that
\begin{align*} 
u^k_l=\nu(q^{k}_l-f_l(z^{k}_l))
\end{align*}
\end{lemma}
\begin{proof}
This follows directly from the optimality condition of $q^k_l$ and Equation \eqref{eq:update u}.
\end{proof}
\begin{lemma}
\label{lemma:u bound}
It holds for every $k\in \mathbb{N}$ and $l=1,\cdots,L-1$ that
\begin{align*}
    \Vert u^{k+1}_l-u^k_l\Vert\leq \nu\Vert q^{k+1}_l-q^{k}_l\Vert+\nu S\Vert z^{k+1}_l-z^{k}_l\Vert
\end{align*}
\begin{proof}
\begin{align*}
    &\Vert u^{k+1}_l-u^k_l\Vert\\&=\Vert \nu(q^{k+1}_l-f_l(z^{k+1}_l))-\nu(q^{k}_l-f_l(z^{k}_l))\Vert \text{(Lemma \ref{lemma: q_opt})}\\&\leq\nu\Vert q^{k+1}_l-q^{k}_l\Vert+\nu\Vert f_l(z^{k+1}_l)-f_l(z^{k}_l)\Vert \text{(triangle inequality)}\\&\leq \nu\Vert q^{k+1}_l-q^{k}_l\Vert+\nu S\Vert z^{k+1}_l-z^{k}_l\Vert \text{(Assumption \ref{ass:lipschitz continuous})}
\end{align*}
\end{proof}
\end{lemma}
\begin{lemma}
\label{lemma:u square bound}
It holds for every $k\in \mathbb{N}$ and $l=1,\cdots,L-1$ that
\begin{align*}
    \Vert u^{k+1}_l-u^k_l\Vert^2_2\leq 2\nu^2(\Vert q^{k+1}_l-q^{k}_l\Vert^2_2+S^2\Vert z^{k+1}_l-z^{k}_l\Vert^2_2)
\end{align*}
\begin{proof}
\begin{align*}
        \Vert u^{k+1}_l-u^k_l\Vert^2_2&=\nu^2\Vert q^{k+1}_l-f_l(z^{k+1}_l)-q^{k}_l+f_l(z^{k}_l)\Vert^2_2 \text{(Lemma \ref{lemma: q_opt})}\\&\leq 2\nu^2(\Vert q^{k+1}_l-q^k_l\Vert^2_2+\Vert f_l(z^{k+1}_l)-f_l(z^{k}_l)\Vert^2_2)\text{(mean inequality)}\\&\leq2\nu^2(\Vert q^{k+1}_l-q^k_l\Vert^2_2+S^2\Vert z^{k+1}_l-z^{k}_l\Vert^2_2)\text{(Assumption \ref{ass:lipschitz continuous})}
\end{align*}

\end{proof}
\end{lemma}
\begin{lemma}
\label{lemma:optimality 1}
For every $k\in \mathbb{N}$, it holds that
\begin{align}  &
  L_\rho(\textbf{p}^{k},\textbf{W}^{k},\textbf{b}^k,\textbf{z}^{k},\textbf{q}^k,\textbf{u}^{k})\!-\! L_\rho(\textbf{p}^{k+1},\textbf{W}^{k},\textbf{b}^k,\textbf{z}^{k},\textbf{q}^k,\textbf{u}^{k}) \geq  \sum\nolimits_{l=2}^{L}(\tau^{k+1}_l/2)\Vert  p^{k+1}_l-p^k_l\Vert^2_2 \label{eq: p optimality} \\&
    L_\rho(\textbf{p}^{k\!+\!1},\textbf{W}^{k},\textbf{b}^k,\textbf{z}^{k},\textbf{q}^k,\textbf{u}^{k})\!-\! L_\rho(\textbf{p}^{k\!+\!1},\textbf{W}^{k\!+\!1},\textbf{b}^{k},\textbf{z}^{k},\textbf{q}^k,\textbf{u}^{k})\! \geq\!  \sum\nolimits_{l\!=\!1}^{L}(\theta^{k\!+\!1}_l/2)\Vert  W^{k\!+\!1}_l\!-\!W^k_l\Vert^2_2 \label{eq: W optimality}
  \\&
  L_\rho(\textbf{p}^{k\!+\!1},\textbf{W}^{k\!+\!1},\textbf{b}^k,\textbf{z}^{k},\textbf{q}^k,\textbf{u}^{k})\!-\! L_\rho(\textbf{p}^{k\!+\!1},\textbf{W}^{k\!+\!1},\textbf{b}^{k\!+\!1},\textbf{z}^{k},\textbf{q}^k,\textbf{u}^{k})\geq(\nu/2) \sum\nolimits_{l\!=\!1}^{L}\Vert b^{k\!+\!1}_l\!-\!b^k_l\Vert^2_2 \label{eq: b optimality}\\& L_\rho(\textbf{p}^{k\!+\!1},\textbf{W}^{k\!+\!1},\textbf{b}^{k\!+\!1},\textbf{z}^{k},\textbf{q}^k,\textbf{u}^{k})\!-\! L_\rho(\textbf{p}^{k\!+\!1},\textbf{W}^{k\!+\!1},\textbf{b}^{k\!+\!1},\textbf{z}^{k\!+\!1},\textbf{q}^k,\textbf{u}^{k})\geq(\nu/2) \sum\nolimits_{l\!=\!1}^{L}\Vert z^{k\!+\!1}_l\!-\!z^k_l\Vert^2_2 \label{eq: z optimality}
\end{align}
\end{lemma}
\begin{proof}
Generally, all inequalities can be obtained by applying optimality conditions of updating $\textbf{p}$, $\textbf{W}$, $\textbf{b}$ and $\textbf{z}$, respectively. We only prove Inequalities \eqref{eq: p optimality},  \eqref{eq: b optimality} and \eqref{eq: z optimality}.This is because Inequality \eqref{eq: W optimality}   follows the same routine of Inequality \eqref{eq: p optimality}.\\
\indent Firstly, we focus on Inequality \eqref{eq: p optimality}. The choice of $\tau^{k+1}_l$ requires
\begin{align}
\phi(p^{k+1}_l,W^k_l,b^k_l,z^k_l,q^k_{l-1},u^k_{l-1})\leq U_l(p^{k+1}_l;\tau^{k+1}_l) \label{eq:stop condition}   
\end{align}
Moreover, the optimality condition of Equation \eqref{eq:update p} leads to
\begin{align}
    \nabla_{p^k_l}\phi(p^k_l,W^k_l,b^k_l,z^k_l,q^k_{l-1},u^k_{l-1})+\tau^{k+1}_l(p^{k+1}_l-p^k_l)=0 \label{eq:p optimality condition}
\end{align}
Therefore
\begin{align*}
    &L_\rho(\textbf{p}^{k},\textbf{W}^{k},\textbf{b}^k,\textbf{z}^{k},\textbf{q}^k,\textbf{u}^{k})\!-\! L_\rho(\textbf{p}^{k+1},\textbf{W}^{k},\textbf{b}^k,\textbf{z}^{k},\textbf{q}^k,\textbf{u}^{k})\\&=\sum\nolimits_{l=2}^L(\phi(p^{k}_l,W^k_l,b^k_l,z^k_l,q^k_{l-1},u^k_{l-1})-\phi(p^{k+1}_l,W^k_l,b^k_l,z^k_l,q^k_{l-1},u^k_{l-1}))\\&\geq\sum\nolimits_{l=2}^L (\phi(p^{k}_l,W^k_l,b^k_l,z^k_l,q^k_{l-1},u^k_{l-1})-U_l(p^{k+1}_l;\tau^{k+1}_l)) \text{(Inequality \eqref{eq:stop condition})}\\&=\sum\nolimits_{l=2}^L(-\nabla_{p^k_l}\phi^T(p^k_l,W^k_l,b^k_l,z^k_l,q^k_{l-1},u^k_{l-1})(p^{k+1}_l-p^k_l)-({\tau}^{k+1}_l/2)\Vert p^{k+1}_{l}-p^k_{l}\Vert^2_2)\\&=\sum\nolimits_{l=2}^L (\tau^{k+1}_l/2)\Vert p^{k+1}_l-p^k_l\Vert^2_2\text{(Equation \eqref{eq:p optimality condition})}
\end{align*}
Next, we prove Inequality \eqref{eq: b optimality}. Because $\nabla_{{b}_1}\phi({p}_1,{W}_1,{b}_1,{z}_1)$ and  $\nabla_{{b}_l}\phi({p}_l,{W}_l,{b}_l,{z}_l,{q}_l,{u}_l)$ are   Lipschitz continuous with coefficient $\nu$. According to Lemma 2.1 in \cite{beck2009fast}, we have
\begin{align}
\nonumber
    \phi({p}^{k\!+\!1}_1,{W}^{k+1}_1,{b}^{k+1}_1,{z}^{k}_1)&\leq\phi({p}^{k\!+\!1}_1,{W}^{k+1}_1,{b}^{k}_1,{z}^{k}_1)+\nabla_{{b}^{k}_1} \phi^T({p}^{k\!+\!1}_1,{W}^{k+1}_1,{b}^{k}_1,{z}^{k}_1)(b^{k+1}_1-b^k_1)\\&+(\nu/2)\Vert b_1^{k+1}-b_1^{k}\Vert^2_2 \label{eq:b1 lipschitz}
    \\\nonumber
    \phi({p}^{k\!+\!1}_l,{W}^{k+1}_l,{b}^{k+1}_l,{z}^{k}_l,{q}^{k}_{l-1},{u}^{k}_{l-1})&\leq\phi({p}^{k\!+\!1}_l,{W}^{k+1}_l,{b}^{k}_l,{z}^{k}_l,{q}^{k}_{l-1},{u}^{k}_{l-1})\\&\nonumber+\nabla_{{b}^{k}_l} \phi^T({p}^{k\!+\!1}_l,{W}^{k+1}_l,{b}^{k}_l,{z}^{k}_l,{q}^{k}_{l-1},{u}^{k}_{l-1})(b^{k+1}_l-b^k_l)\\&+(\nu/2)\Vert b_l^{k+1}-b_l^{k}\Vert^2_2
    \label{eq:b lipschitz}
\end{align}
Moreover, the optimality condition of Equation \eqref{eq:update b} leads to
\begin{align}
    \nabla_{b^k_1}\phi(p^k_1,W^k_1,b^k_1,z^k_1)+\nu(b^{k+1}_1-b^k_1)=0 \label{eq:b1 optimality condition}\\
    \nabla_{b^k_l}\phi(p^k_l,W^k_l,b^k_l,z^k_l,q^k_{l-1},u^k_{l-1})+\nu(b^{k+1}_l-b^k_l)=0 \label{eq:b optimality condition}
\end{align}
Therefore, we have
\begin{align*}
    &L_\rho(\textbf{p}^{k\!+\!1},\textbf{W}^{k+1},\textbf{b}^k,\textbf{z}^{k},\textbf{q}^k,\textbf{u}^{k})\!-\! L_\rho(\textbf{p}^{k\!+\!1},\textbf{W}^{k\!+\!1},\textbf{b}^{k+1},\textbf{z}^{k},\textbf{q}^k,\textbf{u}^{k})\\&=\phi({p}^{k\!+\!1}_1,{W}^{k+1}_1,{b}^k_1,{z}^{k}_1)-\phi({p}^{k\!+\!1}_1,{W}^{k+1}_1,{b}^{k+1}_1,{z}^{k}_1)\\&+\sum\nolimits_{l=2}^L (\phi({p}^{k\!+\!1}_l,{W}^{k+1}_l,{b}^k_l,{z}^{k}_l,q^k_{l-1},u^k_{l-1})-\phi({p}^{k\!+\!1}_l,{W}^{k+1}_l,{b}^{k+1}_l,{z}^{k}_l,q^k_{l-1},u^k_{l-1}))\\&\geq -\nabla_{{b}^{k}_1} \phi^T({p}^{k\!+\!1}_1,{W}^{k+1}_1,{b}^{k}_1,{z}^{k}_1)(b^{k+1}_1-b^k_1)-(\nu/2)\Vert b_1^{k+1}-b_1^{k}\Vert^2_2\\&+\sum\nolimits_{l=2}^L (-\nabla_{{b}^{k}_l} \phi^T({p}^{k\!+\!1}_l,{W}^{k+1}_l,{b}^{k}_l,{z}^{k}_l,{q}^{k}_{l-1},{u}^{k}_{l-1})(b^{k+1}_l-b^k_l)-(\nu/2)\Vert b_l^{k+1}-b_l^{k}\Vert^2_2)\\&\text{(Inequalities \eqref{eq:b1 lipschitz} and \eqref{eq:b lipschitz})}\\&=(\nu/2)\sum\nolimits_{l=1}^L\Vert b^{k+1}_l-b^k_l\Vert^2_2 \text{(Equations \eqref{eq:b1 optimality condition} and \eqref{eq:b optimality condition})}
\end{align*}
Finally, we prove Inequality \eqref{eq: z optimality}. Because $z^{k+1}_l$ minimizes Equation \eqref{eq:update z}, we have
\begin{align}
    \nonumber &(\nu/2)\Vert z^{k+1}_l-W^{k+1}_lp^{k+1}_l-b^{k+1}_l\Vert^2_2+(\nu/2)\Vert q^k_l-f_l(z^{k+1}_l)\Vert^2_2+(\nu/2)\Vert z^{k+1}_l-z^k_l\Vert^2_2\\&\leq (\nu/2)\Vert z^{k}_l-W^{k+1}_l p^{k+1}_l-b^{k+1}_l\Vert^2_2+(\nu/2)\Vert q^k_l-f_l(z^{k}_l)\Vert^2_2 \label{eq:z minimize}
\end{align}
And
\begin{align}
    \nonumber &R(z^k_L;y)\!+(\nu/2)\!\Vert z^k_L\!-\!W^{k\!+\!1}_Lp^{k\!+\!1}_L\!-\!b^{k\!+\!1}_L\Vert^2_2-R(z^{k+1}_L;y)\!-(\nu/2)\!\Vert z^{k+1}_L\!-\!W^{k\!+\!1}_Lp^{k\!+\!1}_L\!-\!b^{k\!+\!1}_L\Vert^2_2\\\nonumber&=R(z^k_L;y)-R(z^{k+1}_L;y)+(\nu/2)\Vert z^{k}_L-z^{k+1}_L\Vert^2_2+\nu(z^{k+1}_L-W^{k+1}_Lp^{k+1}_L-b^{k+1}_L)^T(z^{k}_L-z_L^{k+1})\\\nonumber&\text{($\Vert a-b\Vert^2_2-\Vert a-c\Vert^2_2=\Vert b-c\Vert^2_2+2(c-a)^T(b-c)$ where $a=W^{k+1}_Lp^{k+1}_L+b^{k+1}_L$, $b=z^k_L$, and $c=z^{k+1}_L$)}\\\nonumber &\geq s^T(z^k_L-z^{k+1}_L)+(\nu/2)\Vert z^{k}_L-z^{k+1}_L\Vert^2_2+\nu(z^{k+1}_L-W^{k+1}_Lp^{k+1}_L-b^{k+1}_L)^T(z^{k}_L-z_L^{k+1})\\&\nonumber\text{($s\in \partial R(z^{k+1}_L;y)$ is a subgradient of $R(z^{k+1}_L;y)$)}\\&=(\nu/2)\Vert z^{k+1}_L-z^{k}_L\Vert^2_2\label{eq:zl minimize}\\\nonumber& \text{($0\in s+\nu(z^{k+1}_L-W^{k+1}_Lp^{k+1}_L-b^{k+1}_L)$ by the optimality condition of Equation \eqref{eq:update zl})}
\end{align}
Therefore
\begin{align*}
 &L_\rho(\textbf{p}^{k\!+\!1},\textbf{W}^{k\!+\!1},\textbf{b}^{k\!+\!1},\textbf{z}^{k},\textbf{q}^k,\textbf{u}^{k})\!-\! L_\rho(\textbf{p}^{k\!+\!1},\textbf{W}^{k\!+\!1},\textbf{b}^{k\!+\!1},\textbf{z}^{k\!+\!1},\textbf{q}^k,\textbf{u}^{k})\\&=\sum\nolimits_{i=1}^{L-1}    ((\nu/2)\Vert z^{k}_l-W^{k+1}_l p^{k+1}_l-b^{k+1}_l\Vert^2_2+(\nu/2)\Vert q^k_l-f_l(z^{k}_l)\Vert^2_2\\&-(\nu/2)\Vert z^{k+1}_l-W^{k+1}_lp^{k+1}_l-b^{k+1}_l\Vert^2_2-(\nu/2)\Vert q^k_l-f_l(z^{k+1}_l)\Vert^2_2)\\&+R(z^k_L;y)\!+(\nu/2)\!\Vert z^k_L\!-\!W^{k\!+\!1}_Lp^{k\!+\!1}_L\!-\!b^{k\!+\!1}_L\Vert^2_2-R(z^{k+1}_L;y)\!-(\nu/2)\!\Vert z^{k+1}_L\!-\!W^{k\!+\!1}_Lp^{k\!+\!1}_L\!-\!b^{k\!+\!1}_L\Vert^2_2\\&\geq(\nu/2)\sum\nolimits_{l=1}^L\Vert z^{k+1}_l-z^{k}_l\Vert^2_2 \text{(Inequalities \eqref{eq:z minimize} and \eqref{eq:zl minimize})}
\end{align*}
\end{proof}
\begin{lemma}
For every $k\in N$, it holds that
\begin{align}
\nonumber &L_\rho(\textbf{p}^{k\!+\!1},\textbf{W}^{k\!+\!1},\textbf{b}^{k\!+\!1},\textbf{z}^{k+1},\textbf{q}^k,\textbf{u}^{k})\!-\! L_\rho(\textbf{p}^{k\!+\!1},\textbf{W}^{k\!+\!1},\textbf{b}^{k\!+\!1},\textbf{z}^{k\!+\!1},\textbf{q}^{k+1},\textbf{u}^{k+1})\\&\geq\sum\nolimits_{l=1}^{L-1}((\rho/2-2\nu^2/\rho-\nu/2)\Vert q^{k+1}_l-q^k_l\Vert^2_2-(2\nu^2 S^2/\rho)\Vert z^{k+1}_l-z^k_l\Vert^2_2)
\label{eq:q optimality}
\end{align}
\end{lemma}
\begin{proof}
\begin{align*}
    &L_\rho(\textbf{p}^{k\!+\!1},\textbf{W}^{k\!+\!1},\textbf{b}^{k\!+\!1},\textbf{z}^{k+1},\textbf{q}^k,\textbf{u}^{k})\!-\! L_\rho(\textbf{p}^{k\!+\!1},\textbf{W}^{k\!+\!1},\textbf{b}^{k\!+\!1},\textbf{z}^{k\!+\!1},\textbf{q}^{k+1},\textbf{u}^{k+1})\\&=\sum\nolimits_{l=1}^{L-1} ((\nu/2)\Vert f_l(z^{k+1}_l)-q^k_l\Vert^2_2-(\nu/2)\Vert f_l(z^{k+1}_l)-q^{k+1}_l\Vert^2_2-(u^{k+1}_l)^T(q^k_l-q^{k+1}_l)\\&+(\rho/2)\Vert q^{k+1}_l-q^{k}_l\Vert^2_2-(1/\rho)\Vert u^{k+1}_l-u^k_l\Vert^2_2)\\&=\sum\nolimits_{l=1}^{L-1} ((\nu/2)\Vert f_l(z^{k+1}_l)-q^k_l\Vert^2_2-(\nu/2)\Vert f_l(z^{k+1}_l)-q^{k+1}_l\Vert^2_2-\nu(q^{k+1}_l-f_l(z^{k+1}_l))^T(q^k_l-q^{k+1}_l)\\&+(\rho/2)\Vert q^{k+1}_l-q^{k}_l\Vert^2_2-(1/\rho)\Vert u^{k+1}_l-u^k_l\Vert^2_2)\text{(Lemma \ref{lemma: q_opt})}\\&\geq\sum\nolimits_{l=1}^{L-1}( -(\nu/2)\Vert q^{k+1}_l-q^k_l\Vert^2_2+(\rho/2)\Vert q^{k+1}_l-q^k_l\Vert^2_2-(1/\rho)\Vert u^{k+1}_l-u^k_l\Vert^2_2)\\&\text{($ -\nu(q_l-f_l(z^{k+1}_l))=-(\nu/2)\nabla_{q_l}\Vert q_l-f_l(z^{k+1}_l)\Vert^2_2$ is lipschitz continuous with regard to $q_l$ and Lemma 2.1 in \cite{beck2009fast})}\\&\geq\sum\nolimits_{l=1}^{L-1}( -(\nu/2)\Vert q^{k+1}_l-q^k_l\Vert^2_2+(\rho/2)\Vert q^{k+1}_l-q^k_l\Vert^2_2-(2\nu^2/\rho)\Vert q^{k+1}_l-q^k_l\Vert^2_2-(2\nu^2 S^2/\rho)\Vert z^{k+1}_l-z^k_l\Vert^2_2)\\&\text{(Lemma \ref{lemma:u square bound})}\\&=\sum\nolimits_{l=1}^{L-1}((\rho/2-2\nu^2/\rho-\nu/2)\Vert q^{k+1}_l-q^k_l\Vert^2_2-(2\nu^2 S^2/\rho)\Vert z^{k+1}_l-z^k_l\Vert^2_2)\end{align*}
\end{proof}
\textbf{Proof of Lemma \ref{lemma:objective decrease}}
\begin{proof}
We sum up Inequalities \eqref{eq: p optimality}, \eqref{eq: W optimality}, \eqref{eq: b optimality}, \eqref{eq: z optimality}, and \eqref{eq:q optimality} to obtain Inequality \eqref{eq:objective decrease}.
\end{proof}
\textbf{Proof of Lemma \ref{lemma:lower bounded}}
\begin{proof}
There exists $\textbf{q}^{'}$ such that $p^k_{l+1}=q^{'}_l$ and 
\begin{align*}
    F(\textbf{p}^k,\textbf{W}^k,\textbf{b}^k,\textbf{z}^k,\textbf{q}^{'})\geq \min\nolimits_{\textbf{p},\textbf{W},\textbf{b},\textbf{z},\textbf{q}}\{F(\textbf{p},\textbf{W},\textbf{b},\textbf{z},\textbf{q})| p_{l+1}=q_l\}>-\infty
\end{align*}
Therefore, we have
\begin{align*}
    &L_\rho(\textbf{p}^k,\textbf{W}^k,\textbf{b}^k,\textbf{z}^k,\textbf{q}^{k},\textbf{u}^k)=F(\textbf{p}^k,\textbf{W}^k,\textbf{b}^k,\textbf{z}^k,\textbf{q}^{k})+\sum\nolimits_{l=1}^L (u^k_l)^T(p^k_{l+1}-q^k_l)+(\rho/2)\Vert p^k_{l+1}-q^k_l \Vert^2_2\\&=R(z^k_L;y)+(\nu/2)(\sum\nolimits_{l=1}^{L}\Vert z^k_l-W^k_lp^k_l-b^k_l\Vert^2_2+\sum\nolimits_{l=1}^{L-1}\Vert q^k_l-f_l(z^k_l)\Vert^2_2)\\&+\sum\nolimits_{l=1}^{L-1}( (u^k_l)^T(p^k_{l+1}-q^k_l)+(\rho/2)\Vert p^k_{l+1}-q^k_l \Vert^2_2)\\&=R(z^k_L;y)+(\nu/2)(\sum\nolimits_{l=1}^{L}\Vert z^k_l-W^k_lp^k_l-b^k_l\Vert^2_2+\sum\nolimits_{l=1}^{L-1}\Vert q^k_l-f_l(z^k_l)\Vert^2_2)\\&+\sum\nolimits_{l=1}^{L-1}( \nu(q^k_l-f_l(z^k_l))^T(q^{'}_l-q^k_l)+(\rho/2)\Vert p^k_{l+1}-q^k_l \Vert^2_2)\\&\text{($p^k_{l+1}=q^{'}_l$ and Lemma \ref{lemma: q_opt})}\\&\geq R(z^k_L;y)+(\nu/2)(\sum\nolimits_{l=1}^{L}\Vert z^k_l-W^k_lp^k_l-b^k_l\Vert^2_2+\sum\nolimits_{l=1}^{L-1}\Vert q^{'}_l-f_l(z^k_l)\Vert^2_2)\\&-\sum\nolimits_{l=1}^{L-1}(\nu/2)\Vert q^{'}_l-q^k_l\Vert^2_2 +\sum\nolimits_{l=1}^{L-1}(\rho/2)\Vert p^k_{l+1}-q^k_l \Vert^2_2)\\&\text{($ \nu(q_l-f_l(z^{k+1}_l))=(\nu/2)\nabla_{q_l}\Vert q_l-f_l(z^{k+1}_l)\Vert^2_2$ is lipschitz continuous with regard to $q_l$ and Lemma 2.1 in \cite{beck2009fast})}\\&= F(\textbf{p}^k,\textbf{W}^k,\textbf{b}^k,\textbf{z}^k,\textbf{q}^{'})+(\rho-\nu)/2\Vert p^k_{l+1}-q^k_l\Vert^2_2> -\infty
\end{align*}
Therefore, $F(\textbf{p}^k,\textbf{W}^k,\textbf{b}^k,\textbf{z}^k,\textbf{q}^{'})$ and $(\rho-\nu)/2\Vert p^k_{l+1}-q^k_l\Vert^2_2$ are upper bounded by $L_\rho(\textbf{p}^k,\textbf{W}^k,\textbf{b}^k,\textbf{z}^k,\textbf{q}^{k},\textbf{u}^k)$ and hence $L_\rho(\textbf{p}^0,\textbf{W}^0,\textbf{b}^0,\textbf{z}^0,\textbf{q}^{0},\textbf{u}^0)$ (Lemma \ref{lemma:objective decrease}).
From Assumption \ref{ass:lipschitz continuous}, $(\textbf{p}^k,\textbf{W}^k,\textbf{b}^k,\textbf{z}^k)$ is bounded. $\textbf{q}^k$ is also bounded because $(\rho-\nu)/2\Vert p^k_{l+1}-q^k_l\Vert^2_2$ is upper bounded. $\textbf{u}^k$ is bounded because of Lemma \ref{lemma: q_opt}.
\end{proof}
\textbf{Proof of Lemma \ref{lemma:subgradient bound}}
\begin{proof} We know that
$\partial L_\rho(\textbf{p}^{k+1},\textbf{W}^{k+1},\textbf{b}^{k+1},\textbf{z}^{k+1},\textbf{q}^{k+1},\textbf{u}^{k+1})=\{\nabla_{\text{p}^{k+1}}L_\rho,\nabla_{\text{W}^{k+1}}L_\rho,\nabla_{\text{b}^{k+1}}L_\rho,\partial_{\text{z}^{k+1}}L_\rho,\nabla_{\text{q}^{k+1}}L_\rho,\nabla_{\text{u}^{k+1}}L_\rho\}$ \cite{wang2019admm}. Specifically, we prove that $\Vert g\Vert$ is upper bounded by the linear combination of $\Vert\textbf{p}^{k+1}-\textbf{p}^{k}\Vert$,$ \Vert\textbf{W}^{k+1}-\textbf{W}^{k}\Vert$, $\Vert\textbf{b}^{k+1}-\textbf{b}^{k}\Vert$, $\Vert\textbf{z}^{k+1}-\textbf{z}^{k}\Vert$, $\Vert\textbf{q}^{k+1}-\textbf{q}^{k}\Vert$, and $\Vert\textbf{u}^{k+1}-\textbf{u}^{k}\Vert$.\\
For $p^{k+1}_l$,
\begin{align*}
    &\nabla_{p^{k+1}_l}
L_\rho(\textbf{p}^{k+1},\textbf{W}^{k+1},\textbf{b}^{k+1},\textbf{z}^{k+1},\textbf{q}^{k+1},\textbf{u}^{k+1})\\&=\nabla_{p^{k+1}_l}\phi({p}^{k+1}_l,{W}^{k+1}_l,{b}^{k+1}_l,{z}^{k+1}_l,q^{k+1}_{l-1},u^{k+1}_{l-1})\\&=\nabla_{p^{k}_l}\phi({p}^{k}_l,{W}^{k}_l,{b}^{k}_l,{z}^{k}_l,q^{k}_{l-1},u^{k}_{l-1})+\tau^{k+1}_l(p^{k+1}_l-p^k_l)-\tau^{k+1}_l(p^{k+1}_l-p^k_l)\\&+\nu(W^{k+1}_l)^T W^{k+1}_l p^{k+1}_l-\nu(W^{k}_l)^T W^{k}_l p^{k}_l+\nu(W^{k+1}_l)^T b^{k+1}_l-\nu(W^{k}_l)^T b^{k}_l-\nu(W^{k+1}_l)^T z^{k+1}_l+\nu(W^{k}_l)^T z^{k}_l\\&+(u^{k+1}_{l-1}-u^{k}_{l-1})+\rho(p^{k+1}_l-p^{k}_l)-\rho(q^{k+1}_{l-1}-q^{k}_{l-1})\\&=-\tau^{k+1}_l(p^{k+1}_l-p^k_l)+\nu(W^{k+1}_l)^T W^{k+1}_l p^{k+1}_l-\nu(W^{k}_l)^T W^{k}_l p^{k}_l+\nu(W^{k+1}_l)^T b^{k+1}_l-\nu(W^{k}_l)^T b^{k}_l\\&-\nu(W^{k+1}_l)^T z^{k+1}_l+\nu(W^{k}_l)^T z^{k}_l+(u^{k+1}_{l-1}-u^{k}_{l-1})+\rho(p^{k+1}_l-p^{k}_l)-\rho(q^{k+1}_{l-1}-q^{k}_{l-1})
\end{align*}
So
\begin{align*}
    &\Vert \nabla_{p^{k+1}_l}
L_\rho(\textbf{p}^{k+1},\textbf{W}^{k+1},\textbf{b}^{k+1},\textbf{z}^{k+1},\textbf{q}^{k+1},\textbf{u}^{k+1})\Vert\\&=\Vert \tau^{k+1}_l(p^{k+1}_l-p^k_l)+\nu(W^{k+1}_l)^T W^{k+1}_l p^{k+1}_l-\nu(W^{k}_l)^T W^{k}_l p^{k}_l+\nu(W^{k+1}_l)^T b^{k+1}_l-\nu(W^{k}_l)^T b^{k}_l\\&-\nu(W^{k+1}_l)^T z^{k+1}_l+\nu(W^{k}_l)^T z^{k}_l+(u^{k+1}_{l-1}-u^{k}_{l-1})+\rho(p^{k+1}_l-p^{k}_l)-\rho(q^{k+1}_{l-1}-q^{k}_{l-1})\Vert\\&\leq \tau^{k+1}_l\Vert p^{k+1}_l-p^k_l\Vert+\nu\Vert(W^{k+1}_l)^T W^{k+1}_l p^{k+1}_l-(W^{k}_l)^T W^{k}_l p^{k}_l\Vert+\nu\Vert(W^{k+1}_l)^T b^{k+1}_l-(W^{k}_l)^T b^{k}_l\Vert\\&+\nu\Vert (W^{k+1}_l)^T z^{k+1}_l-(W^{k}_l)^T z^{k}_l\Vert+\Vert u^{k+1}_{l-1}-u^{k}_{l-1}\Vert+\rho\Vert p^{k+1}_l-p^{k}_l\Vert+\rho\Vert q^{k+1}_{l-1}-q^{k}_{l-1}\Vert \ \text{(triangle inequality)}\\&=\tau^{k+1}_l\Vert p^{k+1}_l-p^k_l\Vert+\nu\Vert(W^{k+1}_l)^T W^{k+1}_l (p^{k+1}_l-p^k_l)+(W^{k+1}_l)^T(W^{k+1}_l-W^{k}_l)p^k_l+(W^{k+1}_l-W^{k}_l)^T W^{k}_l p^{k}_l\Vert\\&+\nu\Vert(W^{k+1}_l)^T (b^{k+1}_l-b^k_l)+(W^{k+1}_l-W^{k}_l)^T b^{k}_l\Vert+\nu\Vert(W^{k+1}_l)^T (z^{k+1}_l-z^k_l)+(W^{k+1}_l-W^{k}_l)^T z^{k}_l\Vert\\&+\Vert u^{k+1}_{l-1}-u^k_{l-1}\Vert+\rho\Vert p^{k+1}_l-p^{k}_l\Vert+\rho\Vert q^{k+1}_{l-1}-q^{k}_{l-1}\Vert\\&\leq \tau^{k+1}_l\Vert p^{k+1}_l-p^k_l\Vert+\nu\Vert W^{k+1}_l\Vert^2\Vert p^{k+1}_l-p^k_l\Vert+\nu\Vert W^{k+1}_l\Vert \Vert W^{k+1}_l-W^{k}_l\Vert \Vert p^k_l\Vert+\nu\Vert W^{k+1}_l-W^{k}_l\Vert\Vert W^{k}_l\Vert\Vert p^{k}_l\Vert\\&+\nu\Vert W^{k+1}_l\Vert \Vert b^{k+1}_l-b^k_l\Vert +\nu\Vert W^{k+1}_l-W^{k}_l\Vert\Vert b^{k}_l\Vert+\nu\Vert W^{k+1}_l \Vert \Vert z^{k+1}_l-z^k_l\Vert+\nu\Vert W^{k+1}_l-W^{k}_l\Vert \Vert z^{k}_l\Vert\\&+\nu(\Vert q^{k+1}_{l-1}-q^{k}_{l-1}\Vert+S\Vert z^{k+1}_{l-1}-z^{k}_{l-1}\Vert)+\rho\Vert p^{k+1}_l-p^{k}_l\Vert+\rho\Vert q^{k+1}_{l-1}-q^{k}_{l-1}\Vert\\&\text{(triangle inequality, Cauthy-Schwartz inequality and Lemma \ref{lemma:u bound})}\\&\leq \tau^{k+1}_l\Vert p^{k+1}_l-p^k_l\Vert+\nu\mathbb{N}^2_\textbf{W}\Vert p^{k+1}_l-p^k_l\Vert+2\nu\mathbb{N}_\textbf{W}\mathbb{N}_\textbf{p} \Vert W^{k+1}_l-W^{k}_l\Vert+\nu\mathbb{N}_\textbf{W} \Vert b^{k+1}_l-b^k_l\Vert +\nu\mathbb{N}_\textbf{b}\Vert W^{k+1}_l-W^{k}_l\Vert\\&+\nu\mathbb{N}_\textbf{W} \Vert z^{k+1}_l-z^k_l\Vert+\nu\mathbb{N}_\textbf{z}\Vert W^{k+1}_l-W^{k}_l\Vert+2\nu^2(\Vert q^{k+1}_{l-1}-q^{k}_{l-1}\Vert^2_2+S^2\Vert z^{k+1}_{l-1}-z^{k}_{l-1}\Vert^2_2)+\rho\Vert p^{k+1}_l-p^{k}_l\Vert+\rho\Vert q^{k+1}_{l-1}-q^{k}_{l-1}\Vert\\&\text{(Lemma  \ref{lemma:lower bounded})}
\end{align*}
For $W^{k+1}_1$,
\begin{align*}
&\nabla_{W^{k+1}_1}
L_\rho(\textbf{p}^{k+1},\textbf{W}^{k+1},\textbf{b}^{k+1},\textbf{z}^{k+1},\textbf{q}^{k+1},\textbf{u}^{k+1})\\&=\nabla_{W^{k+1}_1} \phi({p}^{k+1}_1,{W}^{k+1}_1,{b}^{k+1}_1,{z}^{k+1}_1)\\&=\nabla_{W^{k}_1} \phi({p}^{k+1}_1,{W}^{k}_1,{b}^{k}_1,{z}^{k}_1)+\theta^{k+1}_1(W^{k+1}_1-W^k_1)+\nu(W^{k+1}_1-W^k_1)p^{k+1}_1(p^{k+1}_1)^T+\nu(b^{k+1}_1-b^k_1)(p^{k+1}_1)^T\\&-\nu(z^{k+1}_1-z^k_1)(p^{k+1}_1)^T-\theta^{k+1}_1(W^{k+1}_1-W^k_1)\\&=\nu(W^{k+1}_1-W^k_1)p^{k+1}_1(p^{k+1}_1)^T+\nu(b^{k+1}_1-b^k_1)(p^{k+1}_1)^T-\nu(z^{k+1}_1-z^k_1)(p^{k+1}_1)^T-\theta^{k+1}_1(W^{k+1}_1-W^k_1)\\ & \text{(The optimality condition of Equation \eqref{eq:update W})}
\end{align*}
So
\begin{align*}
    &\Vert\nabla_{W^{k+1}_1}
L_\rho(\textbf{p}^{k+1},\textbf{W}^{k+1},\textbf{b}^{k+1},\textbf{z}^{k+1},\textbf{q}^{k+1},\textbf{u}^{k+1})\Vert\\&=\Vert \nu(W^{k+1}_1-W^k_1)p^{k+1}_1(p^{k+1}_1)^T+\nu(b^{k+1}_1-b^k_1)(p^{k+1}_1)^T-\nu(z^{k+1}_1-z^k_1)(p^{k+1}_1)^T-\theta^{k+1}_1(W^{k+1}_1-W^k_1)\Vert\\&\leq \nu\Vert W^{k+1}_1-W^k_1\Vert\Vert p^{k+1}_1\Vert^2+\nu\Vert b^{k+1}_1-b^k_1\Vert\Vert p^{k+1}_1\Vert+\nu\Vert z^{k+1}_1-z^k_1\Vert\Vert p^{k+1}_1\Vert+\theta^{k+1}_1\Vert W^{k+1}_1-W^k_1\Vert\\&\text{(triangle inequality and Cauthy-Schwartz inequality)}\\&\leq \nu\Vert W^{k+1}_1-W^k_1\Vert \mathbb{N}_\textbf{p}^2+\nu\Vert b^{k+1}_1-b^k_1\Vert\mathbb{N}_\textbf{p}+\nu\Vert z^{k+1}_1-z^k_1\Vert\mathbb{N}_\textbf{p}+\theta^{k+1}_1\Vert W^{k+1}_1-W^k_1\Vert\ \text{(Theorem \ref{theo:convergent variable})}
\end{align*}
For $W^{k+1}_l(1<l\leq L)$,
\begin{align*}
&\nabla_{W^{k+1}_l}
L_\rho(\textbf{p}^{k+1},\textbf{W}^{k+1},\textbf{b}^{k+1},\textbf{z}^{k+1},\textbf{q}^{k+1},\textbf{u}^{k+1})\\&=\nabla_{W^{k+1}_l} \phi({p}^{k+1}_l,{W}^{k+1}_l,{b}^{k+1}_l,{z}^{k+1}_l,p^{k+1}_{l-1},u^{k+1}_{l-1})\\&=\nabla_{W^{k}_l} \phi({p}^{k+1}_l,{W}^{k}_l,{b}^{k}_l,{z}^{k}_l,p^{k}_{l-1},u^{k}_{l-1})+\theta^{k+1}_l(W^{k+1}_l-W^k_1)+\nu(W^{k+1}_l-W^k_l)p^{k+1}_l(p^{k+1}_l)^T\\&+\nu(b^{k+1}_l-b^k_l)(p^{k+1}_l)^T-\nu(z^{k+1}_l-z^k_l)(p^{k+1}_l)^T-\theta^{k+1}_l(W^{k+1}_l-W^k_l)\\&=\nu(W^{k+1}_l-W^k_l)p^{k+1}_l(p^{k+1}_l)^T+\nu(b^{k+1}_l-b^k_l)(p^{k+1}_l)^T-\nu(z^{k+1}_l-z^k_l)(p^{k+1}_l)^T-\theta^{k+1}_l(W^{k+1}_l-W^k_l)\\ & \text{(The optimality condition of Equation \eqref{eq:update W})}
\end{align*}
So
\begin{align*}
    &\Vert\nabla_{W^{k+1}_l}
L_\rho(\textbf{p}^{k+1},\textbf{W}^{k+1},\textbf{b}^{k+1},\textbf{z}^{k+1},\textbf{q}^{k+1},\textbf{u}^{k+1})\Vert\\&=\Vert \nu(W^{k+1}_l-W^k_l)p^{k+1}_l(p^{k+1}_l)^T+\nu(b^{k+1}_l-b^k_l)(p^{k+1}_l)^T-\nu(z^{k+1}_l-z^k_l)(p^{k+1}_l)^T-\theta^{k+1}_l(W^{k+1}_l-W^k_l)\Vert\\&\leq \nu\Vert W^{k+1}_l-W^k_l\Vert\Vert p^{k+1}_l\Vert^2+\nu\Vert b^{k+1}_l-b^k_l\Vert\Vert p^{k+1}_l\Vert+\nu\Vert z^{k+1}_l-z^k_l\Vert\Vert p^{k+1}_l\Vert+\theta^{k+1}_l\Vert W^{k+1}_l-W^k_l\Vert\\&\text{(triangle inequality and Cauthy-Schwartz inequality)}\\&\leq \nu\Vert W^{k+1}_l-W^k_l\Vert \mathbb{N}_\textbf{p}^2+\nu\Vert b^{k+1}_l-b^k_l\Vert\mathbb{N}_\textbf{p}+\nu\Vert z^{k+1}_l-z^k_l\Vert\mathbb{N}_\textbf{p}+\theta^{k+1}_l\Vert W^{k+1}_l-W^k_l\Vert\ \text{(Theorem \ref{theo:convergent variable})}
\end{align*}
For $b^{k+1}_1$,
\begin{align*}
 &\nabla_{b^{k+1}_1}
L_\rho(\textbf{p}^{k+1},\textbf{W}^{k+1},\textbf{b}^{k+1},\textbf{z}^{k+1},\textbf{q}^{k+1},\textbf{u}^{k+1})\\&= \nabla_{b^{k+1}_1}
\phi(p^{k+1}_1,W^{k+1}_1,b^{k+1}_1,z^{k+1}_1)\\&= \nabla_{b^{k}_1}
\phi(p^{k+1}_1,W^{k+1}_1,b^{k}_1,z^{k}_1)+\nu(b^{k+1}_1-b^k_1)+\nu(z^k_1-z^{k+1}_1)\\&=\nu(z^k_1-z^{k+1}_1) \ \text{(The optimality condition of Equation \eqref{eq:update b})}
\end{align*}
So $\Vert \nabla_{b^{k+1}_1}
L_\rho(\textbf{p}^{k+1},\textbf{W}^{k+1},\textbf{b}^{k+1},\textbf{z}^{k+1},\textbf{q}^{k+1},\textbf{u}^{k+1})\Vert=\nu\Vert z^{k+1}_1-z^k_1\Vert$.\\
For $b^{k+1}_l(1<l\leq L)$,
\begin{align*}
 &\nabla_{b^{k+1}_l}
L_\rho(\textbf{p}^{k+1},\textbf{W}^{k+1},\textbf{b}^{k+1},\textbf{z}^{k+1},\textbf{q}^{k+1},\textbf{u}^{k+1})\\&= \nabla_{b^{k+1}_l}
\phi(p^{k+1}_l,W^{k+1}_l,b^{k+1}_l,z^{k+1}_l,q^{k}_{l-1},u^{k}_{l-1})\\&= \nabla_{b^{k}_l}
\phi(p^{k+1}_l,W^{k+1}_l,b^{k}_l,z^{k}_l,q^{k}_{l-1},u^{k}_{l-1})+\nu(b^{k+1}_l-b^k_l)+\nu(z^k_l-z^{k+1}_l)\\&=\nu(z^k_l-z^{k+1}_l) \ \text{(The optimality condition of Equation \eqref{eq:update b})}
\end{align*}
So $\Vert \nabla_{b^{k+1}_l}
L_\rho(\textbf{p}^{k+1},\textbf{W}^{k+1},\textbf{b}^{k+1},\textbf{z}^{k+1},\textbf{q}^{k+1},\textbf{u}^{k+1})\Vert=\nu\Vert z^{k+1}_l-z^k_l\Vert$.\\
For $z^{k+1}_l(l<L)$, 
\begin{align*}
&\partial_{z^{k+1}_l}
L_\rho(\textbf{p}^{k+1},\textbf{W}^{k+1},\textbf{b}^{k+1},\textbf{z}^{k+1},\textbf{q}^{k+1},\textbf{u}^{k+1})\\&=\partial_{z^{k+1}_l}
L_\rho(\textbf{p}^{k+1},\textbf{W}^{k+1},\textbf{b}^{k+1},\textbf{z}^{k+1},\textbf{q}^{k},\textbf{u}^{k})+\nu(z_l^{k+1}-z_l^{k})-\nu(z_l^{k+1}-z_l^{k})-\nu\partial f_l(z^{k+1}_l)\circ(q^{k+1}_l-q^k_l)\text{($\circ$ is Hadamard product)}\\&=-\nu(z_l^{k+1}-z_l^{k})-\nu\partial f_l(z^{k+1}_l)\circ(q^{k+1}_l-q^k_l) \ \text{($ 0\in \partial_{z^{k+1}_l}
L_\rho(\textbf{p}^{k+1},\textbf{W}^{k+1},\textbf{b}^{k+1},\textbf{z}^{k+1},\textbf{q}^{k},\textbf{u}^{k})+\nu(z_l^{k+1}-z_l^{k})$)}
\end{align*}
So
\begin{align*}
   &\Vert\partial_{z^{k+1}_l}
L_\rho(\textbf{p}^{k+1},\textbf{W}^{k+1},\textbf{b}^{k+1},\textbf{z}^{k+1},\textbf{q}^{k+1},\textbf{u}^{k+1})\Vert\\&= \Vert -\nu(z_l^{k+1}-z_l^{k})-\nu\partial f_l(z^{k+1}_l)\circ(q^{k+1}_l-q^k_l)\Vert\\&\leq \nu\Vert z^{k+1}_l-z^k_l\Vert+\nu\Vert \partial f_l(z^{k+1}_l)\Vert\Vert q_l^{k+1}-q_l^k\Vert \text{(Cauthy-Schwartz inequality and triangle inequality)}\\&\leq \nu\Vert z^{k+1}_l-z^k_l\Vert+\nu M\Vert q_l^{k+1}-q_l^k\Vert (\text{$\Vert \partial f_l(z^{k+1}_l)\Vert\leq M$})
\end{align*}
For $z^{k+1}_L$,
$\partial_{z^{k+1}_L}
L_\rho(\textbf{p}^{k+1},\textbf{W}^{k+1},\textbf{b}^{k+1},\textbf{z}^{k+1},\textbf{q}^{k+1},\textbf{u}^{k+1})=0$ by the optimality condition of Equation \eqref{eq:update zl}.\\
For $q^{k+1}_l$,
\begin{align*}
&\nabla_{q^{k+1}_l}
L_\rho(\textbf{p}^{k+1},\textbf{W}^{k+1},\textbf{b}^{k+1},\textbf{z}^{k+1},\textbf{q}^{k+1},\textbf{u}^{k+1})\\&=\nabla_{q^{k+1}_l}
L_\rho(\textbf{p}^{k+1},\textbf{W}^{k+1},\textbf{b}^{k+1},\textbf{z}^{k+1},\textbf{q}^{k+1},\textbf{u}^{k})+u^{k+1}_l-u^{k}_l\\&=u^{k+1}_l-u^{k}_l \ \text{($\nabla_{q^{k+1}_l}
L_\rho(\textbf{p}^{k+1},\textbf{W}^{k+1},\textbf{b}^{k+1},\textbf{z}^{k+1},\textbf{q}^{k+1},\textbf{u}^{k})=0$ by the optiamlity condition of Equation \eqref{eq:update q})} 
\end{align*}
So $\Vert\nabla_{q^{k+1}_l}
L_\rho(\textbf{p}^{k+1},\textbf{W}^{k+1},\textbf{b}^{k+1},\textbf{z}^{k+1},\textbf{q}^{k+1},\textbf{u}^{k+1})\Vert=\Vert u^{k+1}_l-u^{k}_l\Vert$.\\
For $u^{k+1}_l$,
\begin{align*}
&\nabla_{u^{k+1}_l} L_\rho(\textbf{p}^{k+1},\textbf{W}^{k+1},\textbf{b}^{k+1},\textbf{z}^{k+1},\textbf{q}^{k+1},\textbf{u}^{k+1})=(p^{k+1}_{l+1}-q^{k+1}_{l})=(u^{k+1}_{l}-u^{k}_{l})/\rho
\end{align*}
So $\Vert\nabla_{u^{k+1}_l} L_\rho(\textbf{p}^{k+1},\textbf{W}^{k+1},\textbf{b}^{k+1},\textbf{z}^{k+1},\textbf{q}^{k+1},\textbf{u}^{k+1})\Vert=\Vert u^{k+1}_{l}-u^{k}_{l}\Vert/\rho$.\\
In summary, we prove that $\nabla_{\text{p}^{k+1}}L_\rho,\nabla_{\text{W}^{k+1}}L_\rho,\nabla_{\text{b}^{k+1}}L_\rho,\partial_{\text{z}^{k+1}}L_\rho,\nabla_{\text{q}^{k+1}}L_\rho,\nabla_{\text{u}^{k+1}}L_\rho$ are upper bounded by the linear combination of $\Vert\textbf{p}^{k+1}-\textbf{p}^{k}\Vert$,$ \Vert\textbf{W}^{k+1}-\textbf{W}^{k}\Vert$, $\Vert\textbf{b}^{k+1}-\textbf{b}^{k}\Vert$, $\Vert\textbf{z}^{k+1}-\textbf{z}^{k}\Vert$, $\Vert\textbf{q}^{k+1}-\textbf{q}^{k}\Vert$, and $\Vert\textbf{u}^{k+1}-\textbf{u}^{k}\Vert$.
\end{proof}
\end{appendix}
\end{document}